\documentclass [letterpaper,11pt]{article}
\pdfoutput=1
\usepackage[T1]{fontenc}
\usepackage[utf8]{inputenc}
\usepackage[english]{babel}
\usepackage[margin=1in]{geometry}

\usepackage{amsmath}
\usepackage{amssymb}
\usepackage{amsfonts}
\usepackage{array}
\usepackage{float}
\usepackage{graphicx}
\usepackage[font=small]{caption}
\usepackage{subcaption}
\usepackage{appendix}
\usepackage{xfrac}
\usepackage{placeins}
\usepackage{bm}
\usepackage{amsthm}
\usepackage[dvipsnames]{xcolor}
\usepackage{comment}
\usepackage{booktabs} 
\usepackage[skins]{tcolorbox}
	\tcbset{commonstyle/.style = {boxrule = 0pt,sharp corners, enhanced jigsaw, nobeforeafter, boxsep = 0pt, left = \fboxsep, right = \fboxsep}}

\usepackage[hidelinks]{hyperref}
\hypersetup{colorlinks=true, citecolor=blue, linkcolor=blue, urlcolor=blue}

  \definecolor{ICES}{RGB}{94, 156, 174}
  \definecolor{ORANGE}{RGB}{191, 87, 0}
  \definecolor{RED}{RGB}{190, 30, 49}
  \definecolor{SUN}{RGB}{227, 81, 51}
  \definecolor{GREEN}{RGB}{0, 171, 86}
  \definecolor{BLUE}{RGB}{11, 78, 179}
  \definecolor{BROWN}{RGB}{122, 80, 40}
  \definecolor{GREY}{RGB}{50, 50, 50}
  \definecolor{TEAL}{RGB}{0, 160, 176}

\newtheorem*{remark}{Remark}
\newtcolorbox{JGcomment}[1][]{commonstyle, #1}

\newcommand{\eq}[1]{(\ref{eq:#1})}
\newcommand{\fig}[1]{\ref{fig:#1}}
\newcommand{\tab}[1]{\ref{tab:#1}}

\newcommand{\bb}[1]{\mathbb{#1}}
\newcommand{\mc}[1]{\mathcal{#1}}
\newcommand{\mf}[1]{\mathfrak{#1}}

\newcommand{\ms}[1]{\mathsf{#1}}
\newcommand{\tc}[1]{\textcolor{#1}}
\newcommand{\tr}[1]{{\check{#1}}}

\newcommand{\p}{\partial}

\newcommand{\ex}{{e_x}}
\newcommand{\ey}{{e_y}}
\newcommand{\ez}{{e_z}}

\newcommand{\LP}{\text{LP}}
\newcommand{\NA}{\text{NA}}

\newcommand{\rcore}{r_\mathrm{core}}
\newcommand{\rclad}{r_\mathrm{clad}}
\newcommand{\ncore}{n_\mathrm{core}}
\newcommand{\nclad}{n_\mathrm{clad}}
\newcommand{\kcore}{k_\mathrm{core}}
\newcommand{\kclad}{k_\mathrm{clad}}
\newcommand{\keff}{k_\mathrm{eff}}
\newcommand{\klp}{k_\mathrm{lp}}

\newcommand{\acoef}{i\omega |J|J^{-1}\eps J^{-T}}
\newcommand{\bcoef}{i \kenv |J|}
\newcommand{\ccoef}{i\omega |J|J^{-1}\mu_0 J^{-T}}

\newcommand{\Gammah}{\Gamma_{\hskip -1pt h}}

\newcommand{\mtmi}{M_{\mathrm{TMI}}}

\newcommand{\lambdas}{\lambda_{\mathrm{s}}}
\newcommand{\lambdap}{\lambda_{\mathrm{p}}}
\newcommand{\omegaj}{\omega_{\mathrm{j}}}

\newcommand{\omegas}{\omega_{\mathrm{s}}}
\newcommand{\omegap}{\omega_{\mathrm{p}}}
\newcommand{\Ej}{E_{\mathrm{j}}}

\newcommand{\Es}{E_{\mathrm{s}}}

\newcommand{\Hj}{H_{\mathrm{j}}}

\newcommand{\Hs}{H_{\mathrm{s}}}

\newcommand{\Pj}{P_{\mathrm{j}}}

\newcommand{\Ps}{P_{\mathrm{s}}}
\newcommand{\Pp}{P_{\mathrm{p}}}
\newcommand{\Ij}{I_{\mathrm{j}}}

\newcommand{\Is}{I_{\mathrm{s}}}
\newcommand{\Ip}{I_{\mathrm{p}}}
\newcommand{\gj}{g_{\mathrm{j}}}

\newcommand{\gs}{g_{\mathrm{s}}}
\newcommand{\gp}{g_{\mathrm{p}}}

\newcommand{\kenv}{\ms{k}}
\newcommand{\kj}{\kenv_{\mathrm{j}}}

\newcommand{\Gammai}{\Gamma_{\hskip -1pt \mathrm{in}}}
\newcommand{\Gammao}{\Gamma_{\hskip -1pt \mathrm{out}}}
\newcommand{\Gammat}{\Gamma_{\hskip -1pt \mathrm{tr}}}
\newcommand{\Omegat}{\Omega_{\mathrm{tr}}}

\newcommand{\Omegac}{\Omega_{\mathrm{c}}}
\newcommand{\Omegapml}{\Omega_{\mathrm{PML}}}

\def\eps{\varepsilon}

\def\grad{\nabla}

\def\curl{\nabla \times}
\def\tcurl{\mathrm{curl}}
\def\div{\nabla \cdot}
\def\tdiv{\mathrm{div}}

\def\hcurl{\nabla_{\hspace{-1pt}h} \times}

\def\hHcurl{H(\mathrm{curl}, \mathcal{T}_h)}

\def\hH1{H^1(\mathcal{T}_h)}
\def\Hcurl{H(\mathrm{curl}, \Omega)}

\def\H1{H^1(\Omega)}
\def\lb{\langle}
\def\rb{\rangle}

\graphicspath{{figures}}

\usepackage[backend=bibtex,
	maxbibnames=99,
	maxcitenames=99,
	giveninits=true]
	{biblatex}
\title{\Large{A Vectorial Envelope Maxwell Formulation for Electromagnetic Waveguides \\ 
with Application to Nonlinear Fiber Optics}}

\author{
Stefan Henneking$^{1,*}$,
Jacob Grosek$^2$, and
Leszek Demkowicz$^1$ \\
{\normalsize
$^1$Oden Institute, The University of Texas at Austin} \\ 
{\normalsize
$^2$Air Force Research Laboratory}
}
{\date{\normalsize \today}}
\begin{document}
\clearpage \maketitle
\thispagestyle{empty}
\renewcommand*{\thefootnote}{\fnsymbol{footnote}}
\footnotetext{$^*$Corresponding author: stefan@oden.utexas.edu}
\renewcommand*{\thefootnote}{\arabic{footnote}}

\addcontentsline{toc}{section}{Abstract}

\begin{abstract}
This article presents an ultraweak discontinuous Petrov--Galerkin (DPG) formulation of the time-harmonic Maxwell equations for the vectorial envelope of the electromagnetic field in a weakly-guiding multi-mode fiber waveguide.
This formulation is derived using an envelope ansatz for the vector-valued electric and magnetic field components, factoring out an oscillatory term of $\exp(-i \kenv z)$ with a user-defined wavenumber $\kenv$, where $z$ is the longitudinal fiber axis and field propagation direction.
The resulting formulation is a modified system of the time-harmonic Maxwell equations for the vectorial envelope of the propagating field.
This envelope is less oscillatory in the $z$-direction than the original field, so that it can be more efficiently discretized and computed, enabling solution of the vectorial DPG Maxwell system for $1000\times$ longer fibers than previously possible.
Different approaches for incorporating a perfectly matched layer for absorbing the outgoing wave modes at the fiber end are derived and compared numerically. 
The resulting formulation is used to solve a 3D Maxwell model of an ytterbium-doped active gain fiber amplifier, coupled with the heat equation for including thermal effects. 
The nonlinear model is then used to simulate thermally-induced transverse mode instability (TMI). 
The numerical experiments demonstrate that it is computationally feasible to perform simulations and analysis of real-length optical fiber laser amplifiers using discretizations of the full vectorial time-harmonic Maxwell equations. 
The approach promises a new high-fidelity methodology for analyzing TMI in high-power fiber laser systems and is extendable to including other nonlinearities. 
\end{abstract}

\vspace{10pt}


\section{Introduction}

\paragraph{Motivation.}
Optical fibers are electromagnetic waveguides that transmit light very efficiently (i.e.~with low losses) over long distances and are useful for a vast number of applications.
Fiber laser amplifiers are optical fiber devices designed for achieving highly coherent light sources with high power outputs.
The combination of high average powers and extremely high beam qualities make fiber lasers an important technology in many industrial, defense, and scientific applications~\cite{jauregui2013fiber}.
However, the efforts of power-scaling beam combinable fiber amplifiers have encountered roadblocks in the form of nonlinear effects~\cite{agrawal, jauregui2013fiber} including the thermally-induced transverse mode instability (TMI).
TMI is characterized by a sudden reduction of quality and stability of the beam emitted by a fiber laser system once a certain power threshold has been reached; this nonlinearity has revealed itself as the one of the strongest limitations for the average power scaling of current fiber laser systems~\cite{jauregui2020tmi}.
Modeling and simulation play a decisive role in understanding the TMI and other detrimental nonlinear effects, finding mitigation strategies, and informing fiber architectures. 

\paragraph{Literature.}
This paper is an extension of the work done in~\cite{nagaraj2018raman, henneking2021fiber, henneking2021phd, henneking2022parallel}. 
Nagaraj et al.~\cite{nagaraj2018raman} developed a 3D vectorial Maxwell model to simulate passive Raman gain amplification in an optical fiber. 
Building on that framework, we added the ability to model the more common active gain amplification through a rare-earth, lanthanide metal dopant in the fiber core region~\cite{henneking2021fiber}; the implementation of this fiber amplifier model was then extended to support large-scale numerical simulations~\cite{henneking2021phd, henneking2022parallel}. 

The model is discretized with the discontinuous Petrov--Galerkin (DPG) finite element method with optimal test functions~\cite{demkowicz2017dpg}. 
As an ultraweak DPG formulation, the model can make use of robust automatic $hp$-adaptive algorithms~\cite{chakraborty2024hp} and advanced solvers for wave propagation~\cite{petrides2021adaptive, badger2023scalable}. 
The computational expense for solving the vectorial Maxwell fiber amplifier model increases slightly more than linearly with the number of wavelengths (which is proportional to the length of the waveguide) due to the effect of numerical pollution~\cite{henneking2021pollution, melenk2023waveguide1, melenk2020maxwell, babuska1997pollution}. 
An efficient parallel implementation of the model, using high-order discretization with fast integration techniques~\cite{mora2019fast, badger2020fast}, is able to simulate thousands of wavelengths in an optical fiber amplifier~\cite{henneking2022parallel}. 
Because of the short wavelength of the light ($\mc{O}(\mu\text{m})$), however, these computations are limited to fiber lengths of only a few centimeters when real-length fiber amplifiers are several meters long.
Given the exceedingly large number of wavelengths within a real-length fiber amplifier, it is computationally intractable to resolve the wavelength scale with a full vectorial Maxwell model. 

A common approach for fiber amplifier modeling is thus to reduce the complexity of the Maxwell equations by introducing additional assumptions or approximations of the physics involved, leading to simplified but computational efficient models such as coupled-mode-theory~\cite{naderi2013tmi, goswami2021fiber} or beam propagation models~\cite{gonthier1991bpm, saitoh2001bpm, ward2013bpm}. 
However, the variety of assumptions made in their derivations may limit their ability to accurately capture some of the nonlinear optical phenomena in fiber amplifiers, which motivates our work on developing Maxwell formulations for optical waveguides that enable higher-fidelity simulations of fiber laser systems. 

\paragraph{Contributions.}
The main contributions of this paper are: 
(1)~derivation of a vectorial envelope DPG formulation of the time-harmonic Maxwell equations that can be efficiently discretized and solved inside a weakly-guiding optical fiber waveguide;
(2)~derivation and analysis of different approaches for implementing absorbing boundaries with a stretched coordinate perfectly matched layer (PML) at the waveguide output; 
(3)~application of the vectorial envelope DPG Maxwell formulation to a fiber laser amplifier model that includes nonlinear gain and thermal effects via coupling to the heat equation; and 
(4)~numerical results of the fiber amplifier model showing the three regimes, \textit{stable}, \textit{transition}, and \textit{chaotic}, of thermally-induced TMI nonlinearity. 
The contributions made in this paper enable efficient computations of the propagating electromagnetic fields inside weakly-guiding multi-mode optical waveguides. 
To the best of our knowledge, this is the first Maxwell model capable of solving for the full vectorial electromagnetic field in a real-length fiber laser amplifier and the first vectorial Maxwell model capable of capturing the onset of the TMI phenomenon.

\section{Background}
\label{sec:background}

\begin{figure}[htb]
	\centering
	\includegraphics[width=\textwidth,trim={0 0 10pt 0},clip]
	{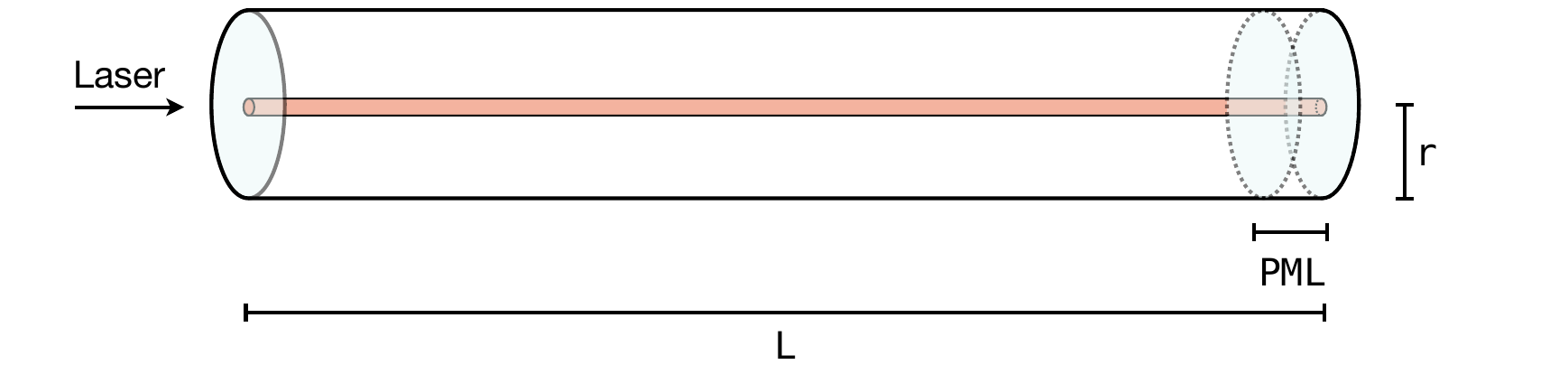}
	\caption{Illustration of a step-index fiber waveguide.}
	\label{fig:fiber-waveguide}
\end{figure}

\paragraph{Maxwell equations in a weakly-guiding step-index fiber.}
A cylindrical step-index glass optical fiber waveguide, as depicted in Figure~\fig{fiber-waveguide}, consists of a cylindrical fiber core region of radius $\rcore$ centered in the fiber, surrounded by a homogeneous cladding region of radius $\rclad$.\footnote{Note that the cladding is itself coated in one or more layers of polymer to protect the glass fiber.}
Suppose that the center of the waveguide is aligned with the z-axis, so that $x, y$ are the transverse directions. 
The material refractive index of the core ($\ncore$) is slightly larger than that of the cladding ($\nclad$), i.e.,
\begin{equation}
\arraycolsep=2pt
	n(r) = 
	\left\{ \begin{array}{ll}
	\ncore & ,\quad r \le \rcore, \\ 
	\nclad < \ncore & ,\quad \rcore < r < \rclad ,
	\end{array} \right.
\end{equation}
where $r = \sqrt{x^2+y^2}$, so as to guide the light in the core region via total internal reflection.
The step-index fiber is called \emph{weakly-guiding} if $(\ncore - \nclad)/\ncore \ll 1$.

Let $\Omega := \Omegat \times (0,L)$ denote the fiber domain, where $\Omegat := \{ (x,y) : x^2 + y^2 < \rclad \}$ is the transverse domain and $L$ is the length of the waveguide.
Let $\Gamma := \p \Omega$; we define the radial boundary, i.e.~the boundary of the transverse domain, $\Gammat := \p \Omegat \times (0,L)$, the fiber input $\Gammai := \Omegat \times \{ 0 \}$, and the fiber output $\Gammao := \Omegat \times \{ L \}$. 
The following discussion partially utilizes the arguments in~\cite{henneking2021phd, jackson, griffiths}.

An optical fiber, primarily comprised of fused silica glass, is a non-magnetic and dielectric medium~\cite{agrawal, snyder1983optical}. In the absence of free charges, the linear time-harmonic Maxwell equations describing the electromagnetic field propagating inside the fiber are given by:
\begin{align}
	\curl E &= -i \omega \mu_0 H , 
	\label{eq:time-harmonic-1} \\
	\curl H &= i \omega \eps E ,
	\label{eq:time-harmonic-2}
\end{align}
where $E$ and $H$ are the electric and magnetic field vectors, $i = \sqrt{-1}$, $\omega$ is the angular frequency, $\mu_0$ is the magnetic permeability in vacuum, and $\eps$ is the electric permittivity.
At the radial boundary ($r = \rclad$), we assume perfect electrical conductor (PEC) boundary conditions (BCs); that is, the tangential electrical field vanishes on $\Gammat$. BCs at the fiber input $\Gammai$ and fiber output $\Gammao$ will be specified later. The first-order system \eq{time-harmonic-1}--\eq{time-harmonic-2} can be reduced to a second-order form to obtain the curl--curl formulation for the electric field:\footnote{One can similarly derive the curl--curl formulation for the magnetic field and then continue with the same arguments to arrive at an analogous eigenvalue problem for the transverse magnetic field components.}
\begin{equation}
	\curl \curl E - \omega^2 \mu_0 \eps E = 0 .
	\label{eq:time-harmonic-curl-curl}
\end{equation}

By applying the vector identity $\curl \curl E = \grad (\div E) - \Delta E$, the curl--curl equation \eq{time-harmonic-curl-curl} can be simplified to a vectorial Helmholtz equation
\emph{under the assumption} that $\grad (\div E) = 0$:
\begin{align}
	\Delta E + \omega^2 \mu_0 \eps E = 0 .
	\label{eq:vectorial-helmholtz}
\end{align}

\begin{remark}{(Vectorial Helmholtz formulations.)}
The assumption $\grad (\div E) = 0$, while common in optical waveguide modeling, is not necessarily a good assumption. In particular, for the fiber amplifier problem (see Section~\ref{sec:amplifier}), when $\eps E = \eps_0 E + P$, where $\eps_0$ is electric permittivity in vacuum, and polarization $P$ is used to model background polarization, gain polarization and thermal polarization, one would have to assume that $\grad (\div P) = 0$, which is not true in general.
\end{remark}

\paragraph{LP modes.}
Assuming a guided wave propagating forward along the $z$-direction, the electric field takes the form:
\begin{align}
	E(x,y,z) = E(x,y) e^{-ikz} ,
	\label{eq:guided-wave-ansatz}
\end{align}
where $k$ is the propagation constant with real part $\mf{Re} \{ k \} > 0$.\footnote{Note that backward traveling waves of the form $E(x,y)e^{+i k z}$ are, in principle, possible as well; they are usually discussed in the context of resonant cavity problems. This implies that our convention for the time-harmonic Maxwell equations~\eqref{eq:time-harmonic-1}--\eqref{eq:time-harmonic-2} was to factor out $e^{+i \omega t}$.}$^{,}$\footnote{When the wavenumber $k$ is determined through the mode eigenproblem~\eqref{eq:mode-eigenproblem}, it is also often referred to as the mode propagation constant.}
Using ansatz \eq{guided-wave-ansatz} for the Helmholtz equation \eq{vectorial-helmholtz}, we obtain a transverse Helmholtz equation for the field envelope:
\begin{align}\label{eq:mode-eigenproblem}
	( \Delta_t + (\omega^2 \mu_0 \eps - k^2) ) E(x,y) = 0,
\end{align}
where $\Delta_t$ is the transverse part of the Laplace operator. Analysis of the corresponding eigenvalue problem, $(\Delta_t + \zeta^2) \Psi = 0$, where $\zeta^2 \equiv \omega^2 \mu_0 \eps - k^2$,
with appropriate boundary conditions, yields a spectrum of positive eigenvalues $\zeta_\lambda^2$, and eigenmodes $\Psi_\lambda, \lambda = 1,2,\ldots$; \emph{guided modes} are those for which the corresponding propagation constant $k_\lambda$ is real-valued; otherwise, the mode is decaying and called evanescent.

We define the cladding and core wavenumbers
\begin{equation}
	\kclad := \frac{\omega}{c} \nclad
	\quad \text{and} \quad
	\kcore := \frac{\omega}{c} \ncore ,
\end{equation}
respectively, where $c = 1/\sqrt{\eps_0 \mu_0}$ is the speed of causality.
Under the weakly-guiding condition, the wave equation may be posed for the transverse electric field components: 
\begin{align}
	\text{Core:} \quad
	\left[\Delta_t + \left( \kcore^2 - k^2 \right) \right]
	\left\{ \begin{array}{c} E_x\\ E_y \end{array} \right\} 
	&= 0,
	\label{eq:helmholtz-core} \\
	\text{Cladding:} \quad
	\left[\Delta_t - \left( k^2 - \kclad^2 \right) \right] 
	\left\{ \begin{array}{c} E_x\\ E_y \end{array} \right\} 
	&= 0 .
	\label{eq:helmholtz-clad}
\end{align}
The corresponding eigenvalue problem with appropriate boundary and core-cladding interface conditions yields transverse \emph{core-guided} modes $\psi_\lambda$ that satisfy 
\begin{equation}
	\kclad < | k_\lambda | < \kcore .
\end{equation}

These modes have two possible linear polarizations in the transverse directions: $\ex$~and~$\ey$. They are therefore called \emph{$\LP$ modes}.
In cylindrical coordinates, they must satisfy the following characteristic equation involving $l$-th order Bessel functions $J_l$ and modified Bessel functions $K_l$ (cf.~\cite[Eqn.~(8.128)]{jackson}):
\begin{equation}
\frac{(\zeta \rcore) J'_l(\zeta \rcore)}{J_l(\zeta \rcore)} = \frac{(\chi \rcore) K'_l(\chi \rcore)}{K_l(\chi \rcore)} , \quad l = 0,1,2, \ldots ,
	\label{eq:characteristic-weakly}
\end{equation}
as well as
\begin{equation}
	(\rcore \zeta)^2 + (\rcore \chi)^2 = \rcore^2 \frac{\omega^2}{c^2} (\ncore^2 - \nclad^2) \equiv V^2 ,
	\label{eq:V}
\end{equation}
where
\begin{equation}
	\zeta^2 = \left( \kcore^2 - k^2 \right), \quad
	\chi^2 = \left( k^2 - \kclad^2 \right) .
\end{equation}
$V$ is called the \emph{normalized frequency} or $V$-number, and $\NA := (\ncore^2 - \nclad^2)^{1/2}$ is the fiber core numerical aperture.

\begin{table}[htb]
		\centering
		\caption[Cutoff frequencies of lowest-order LP modes]{Cutoff frequencies of lowest-order LP modes in a weakly-guiding step-index fiber. The fundamental mode ($\LP_{01}$) has no cutoff and can propagate at any frequency.}
		\label{tab:LP-cutoff}
		\begin{tabular}{@{}rcccccc@{}}
			\toprule
			Guided mode & 
			$\LP_{01}$ & $\LP_{11}$ & $\LP_{21}, \LP_{02}$ & $\LP_{31}$ & $\LP_{12}$ & $\cdots$ \\
			\midrule
			Cutoff frequency $V_c$ & - & 2.405 & 3.832 & 5.136 & 5.520 & $\cdots$ \\
			\bottomrule
		\end{tabular}
\end{table}

\begin{figure}[htb]
	\centering
	\begin{subfigure}[b]{0.87\columnwidth}
		\includegraphics[width=0.32\textwidth]{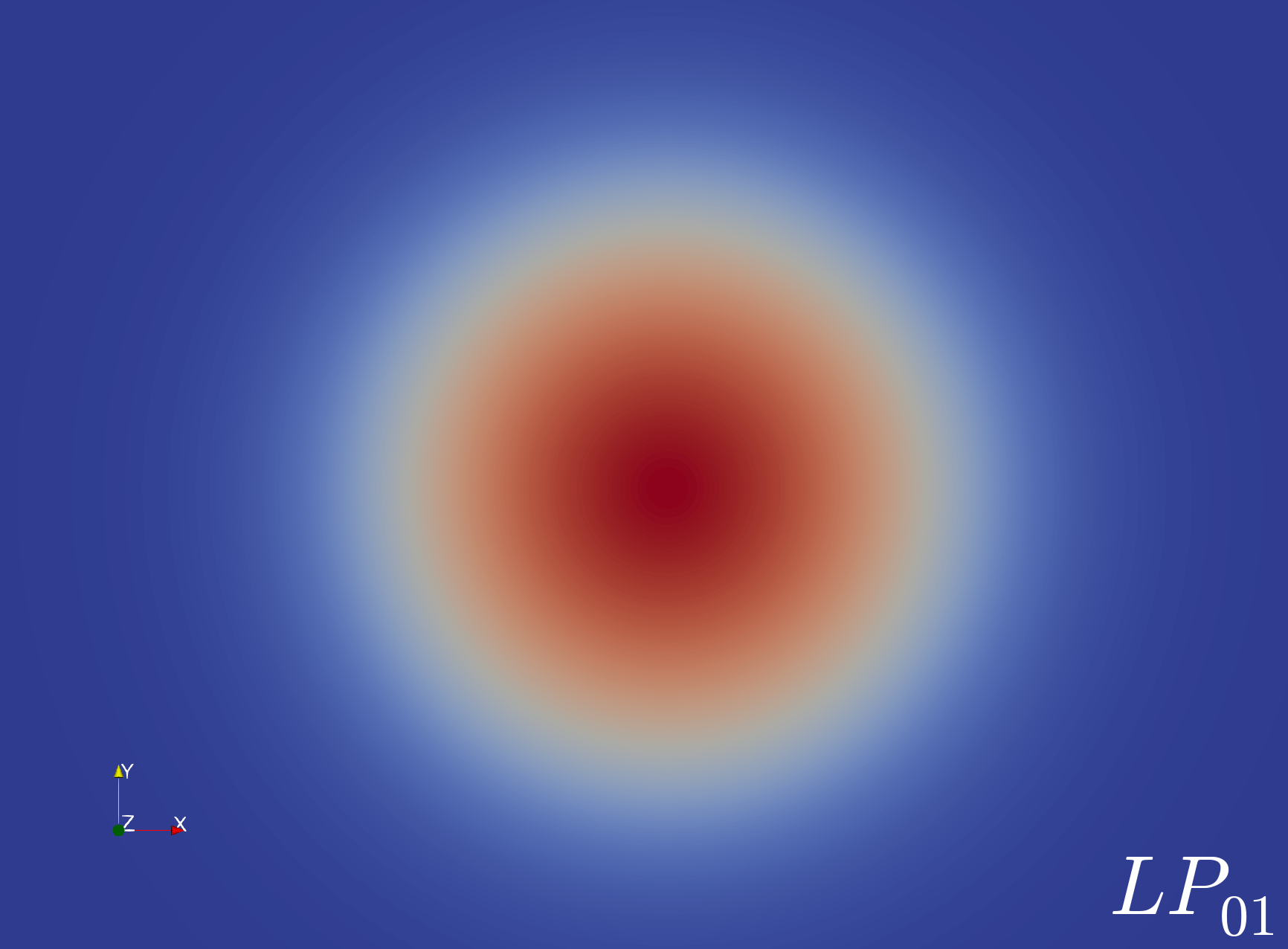}
		\includegraphics[width=0.32\textwidth]{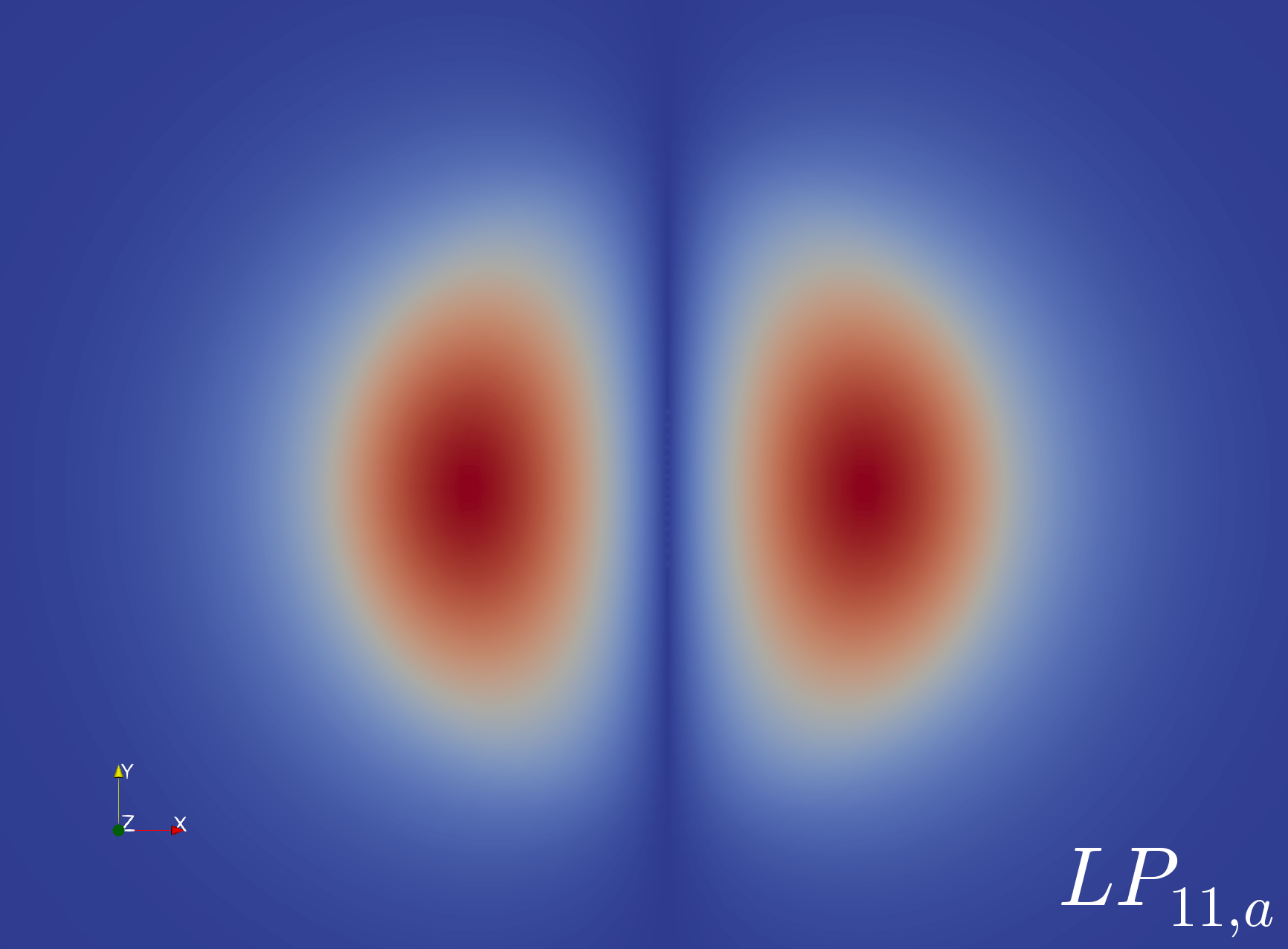}
		\includegraphics[width=0.32\textwidth]{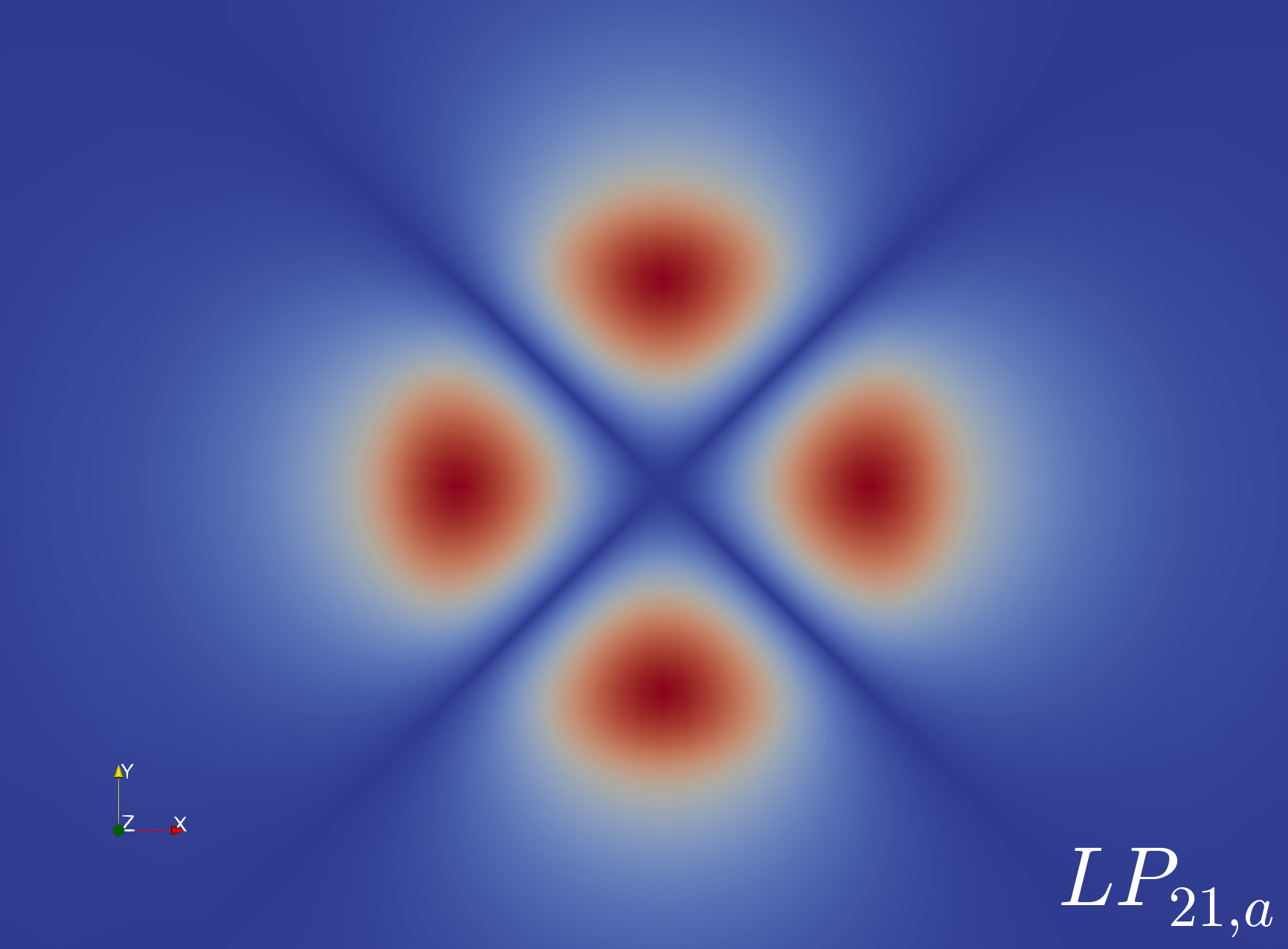}\\
		\includegraphics[width=0.32\textwidth]{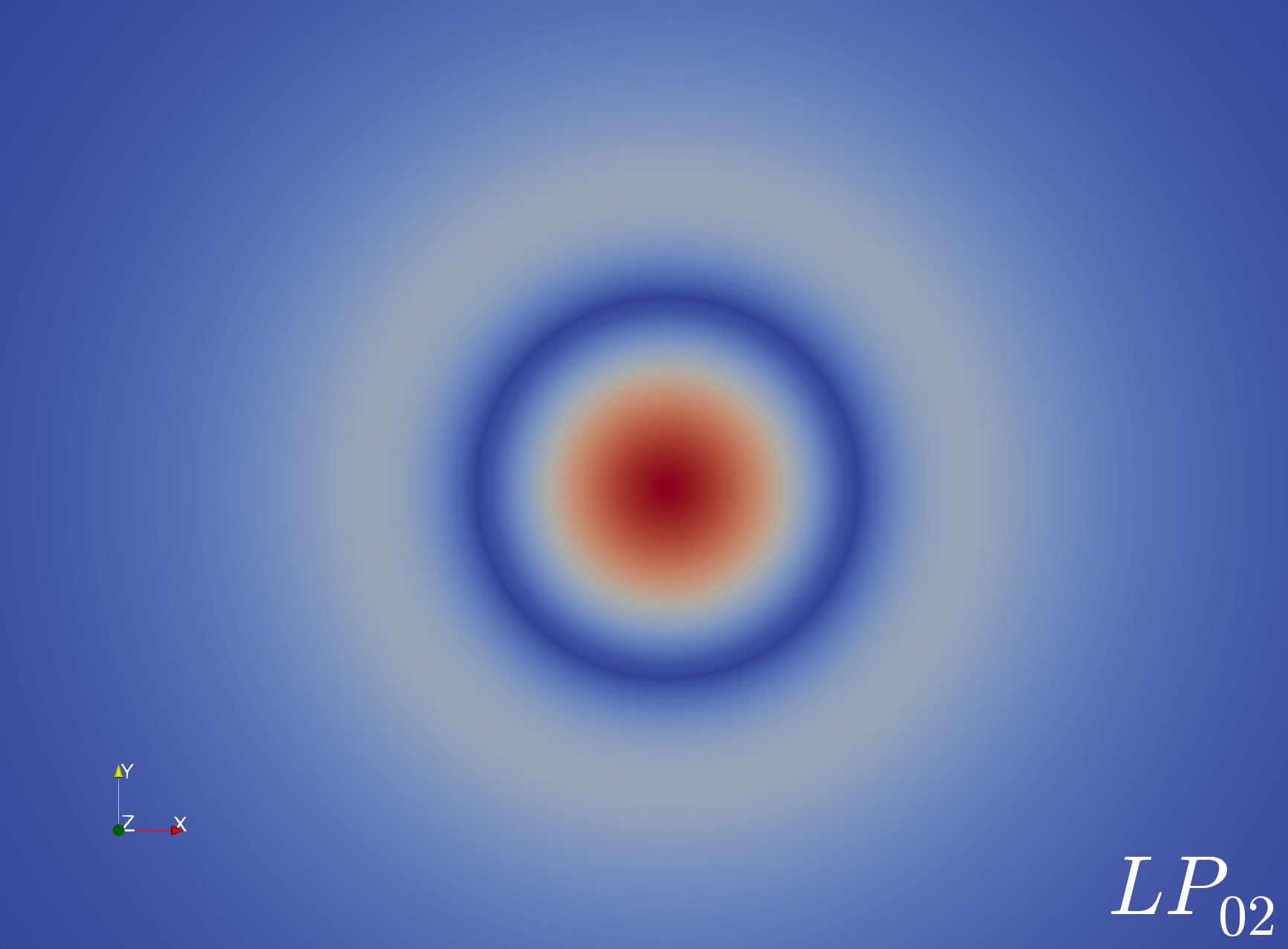}
		\includegraphics[width=0.32\textwidth]{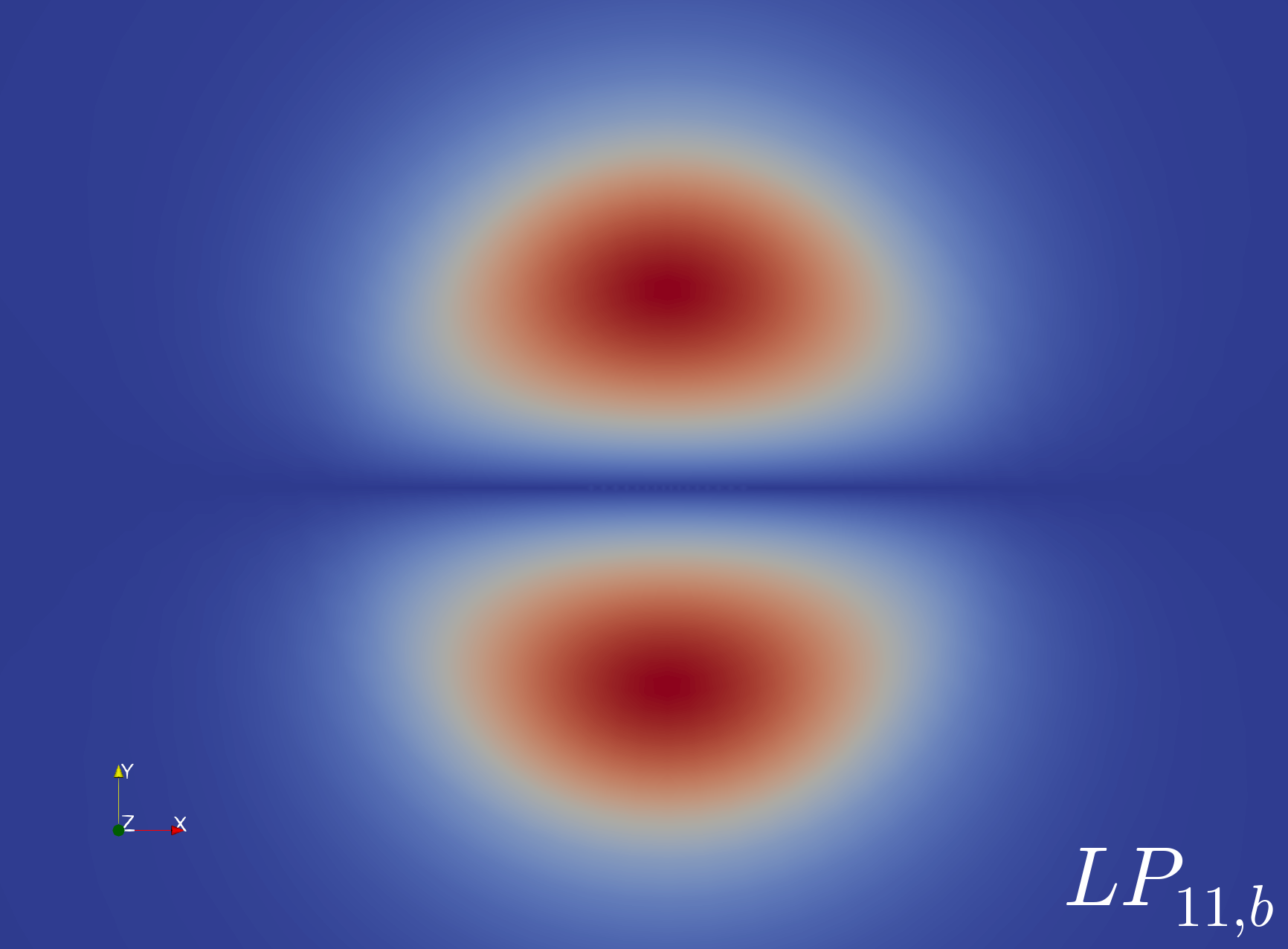}
		\includegraphics[width=0.32\textwidth]{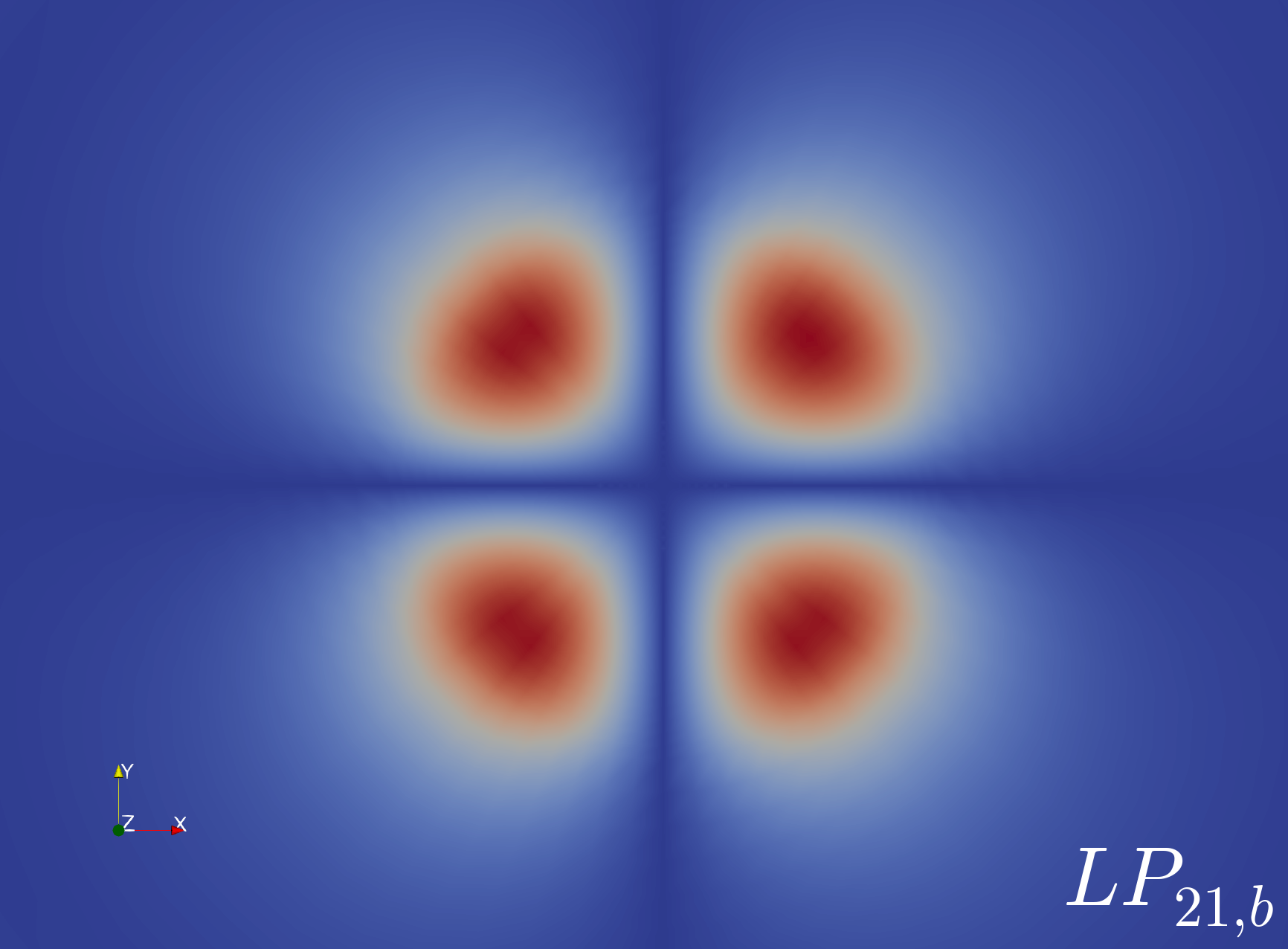}
	\end{subfigure}
	\begin{subfigure}[b]{0.12\columnwidth}
		\includegraphics[width=\textwidth]{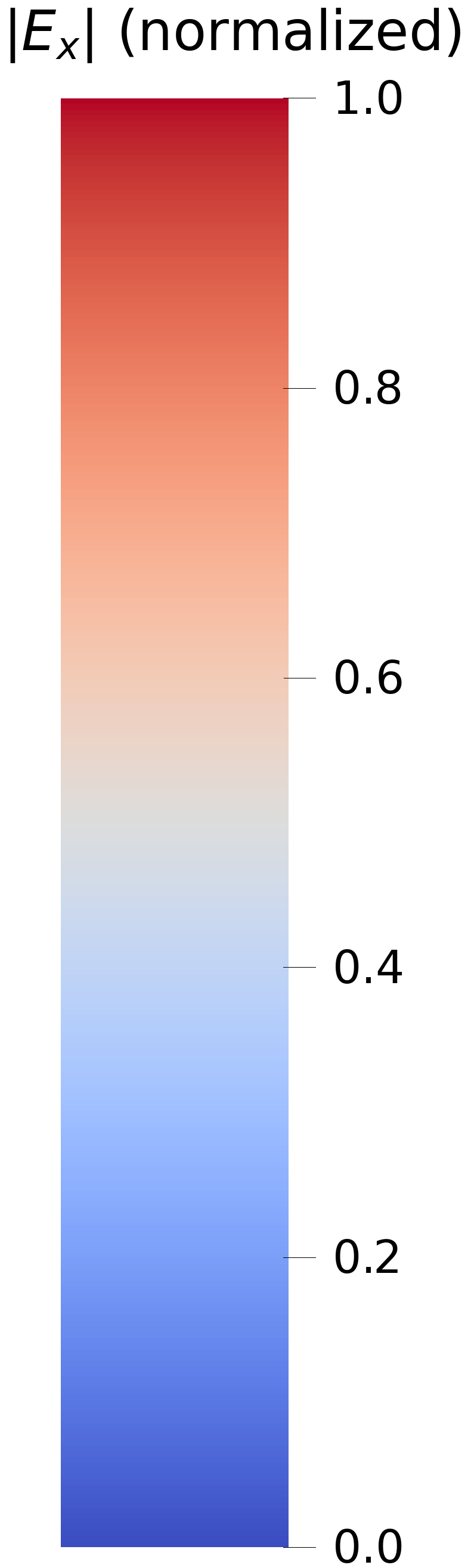}
	\end{subfigure}
	\caption{Electric field magnitude of the lowest-order transverse core-guided LP modes in a weakly-guiding step-index fiber. The $\LP_{01}$ mode is the fundamental mode (FM) and has no cutoff frequency. All other LP modes are higher-order modes (HOMs) that can only propagate above their respective cutoff frequency. Asymmetric HOMs such as $\LP_{11}$ or $\LP_{21}$ have multiple rotations. In a symmetric (i.e.~not birefringent) step-index fiber, these rotated modes all have the same cutoff frequency.}
	\label{fig:LP-modes}
\end{figure}

\begin{table}[htb]
	\centering
	\caption{Step-index fiber parameters}
	\label{tab:fiber-parameters}
	\begin{tabular}{@{}llll@{}}
		\toprule
		Symbol & Description & Value & Unit \\
		\midrule
		$r_{\text{core}}$ & Core radius & $12.7$ & $\mu$m \\
		$r_{\text{clad}}$ & Cladding radius & $127$ & $\mu$m \\
		$n_{\text{core}}$ & Refractive index in fiber core & $1.4512$ & - \\
		$n_{\text{clad}}$ & Refractive index in fiber cladding & $1.4500$ & - \\
		NA & Core numerical aperture & $0.059$ & - \\
		$\lambda$ & Wavelength & 1\,064 & nm \\
		$V(\omega)$ & Normalized frequency & $4.43$ & - \\
		\bottomrule
	\end{tabular}
\end{table}

Given any suitable fiber parameters $\ncore,\nclad,\rcore,\rclad$, and frequency $\omega$, for every $l=0,1,2,\ldots$, there are infinitely many solutions $\zeta, \chi$ that satisfy (\ref{eq:characteristic-weakly}) but only finitely many of these may satisfy (\ref{eq:V}) as well. These are denoted $\zeta_{lp}$, $\chi_{lp}$, $p=1,2,\ldots,N$. We find that only for $l=0$, there exists a solution for any $V > 0$. This \emph{fundamental mode} (FM) is the $\LP_{01}$ mode, and it does not have a cutoff frequency (i.e.~a frequency $V_c$ below which the mode cannot propagate).
All other LP modes are \emph{higher-order modes} (HOMs) that can only propagate if the $V$-number is larger than their respective cutoff frequency~\cite[{\S}12-9]{snyder1983optical}. These cutoff frequencies can be calculated for each mode and are given in Table~\ref{tab:LP-cutoff} for the lowest-order LP modes. 

Table~\ref{tab:fiber-parameters} shows the fiber parameters used for the numerical experiments throughout this paper. With a $V$-number of $4.43$, only four core-guided modes can propagate in the fiber: $\LP_{01}$, $\LP_{11}$, $\LP_{21}$, and $\LP_{02}$.
Figure~\fig{LP-modes} depicts the transverse profile of the (normalized) electric field magnitude of these LP modes.

\paragraph{Definitions and notation.}
For a complex-valued quantity $a \in \bb C$, $\bar a$ denotes the complex conjugate of $a$. Let $(\cdot, \cdot)$ and $\| \cdot \|$ respectively denote the standard $L^2$ inner product (antilinear in the second argument) and associated norm in the Hilbert space $L^2(\Omega)$.
The notation $y \in L^2(\Omega)$ may refer to a function $y$ with a single scalar-valued component in $L^2(\Omega)$ or, if $y$ has multiple components $y = (y_1, y_2, \ldots, y_d)^T$, each component $y_i \in L^2(\Omega), i = 1, 2, \ldots, d$.
We define the space
\begin{equation}
	\Hcurl := \left\{ q : \Omega \rightarrow \bb{C}^3 : q \in L^2(\Omega), \curl q \in L^2(\Omega) \right\} .
\end{equation}

Let $\mc{T}_h$ denote a suitable finite element triangulation of $\Omega$ with mesh skeleton $\Gammah$. The use of discontinuous test functions necessitates an element-local representation of the $\Hcurl$ space:
\begin{equation}
	\hHcurl := \left\{ q : \Omega \rightarrow \bb{C}^3 : q|_K \in H(\tcurl, K)\ \forall K \in \mc{T}_h \right\} \supset \Hcurl .
\end{equation}
Analogously, the notation $\hcurl$ refers to an element-wise interpretation of the operator $\curl$.

Lastly, we must introduce notation related to unknowns defined on the mesh skeleton $\Gammah$. In particular, the trace space on the mesh skeleton is defined through element-wise traces of globally conforming functions:
\begin{equation}
	H^{-1/2}(\tcurl,\Gamma_h) := \gamma_t \Hcurl ,
\end{equation}
where, denoting the outward unit normal by $n$ and the tangential component $q_t := -n \times (n \times q)$,
\begin{equation}
	\gamma_t q := 
	\prod_{K \in \mc{T}_h} \gamma_t^K (q|_K) =
	\prod_{K \in \mc{T}_h} q_t|_{\p K} ,
	\quad q \in \hHcurl .
\end{equation}
The mesh skeleton term $\lb \cdot, \cdot \rb_{\Gammah}$ can then be interpreted as a sum of element-wise duality pairings on element boundaries between the tangential trace $q_t|_{\p K} \in H^{-1/2}(\tcurl,\p K)$ and the rotated tangential trace $n \times \tr{E}_t|_{\p K} \in H^{-1/2}(\tdiv,\p K)$:
\begin{equation}
	\lb n \times \tr E, q \rb_{\Gammah} :=
	\sum_{K \in \mc{T}_h} \lb n \times \tr E_t, q_t
	\rb_{\p K} \, ,
\end{equation}
where $\tr{E} \in H^{-1/2}(\tcurl,\Gamma_h)$ and $q \in \hHcurl$; see \cite[Lem.~2.2]{carstensen2016breaking} and \cite{demkowicz2023fem} for details.\footnote{Instead of the conventional ``hat''-notation for traces (such as $\hat{E}$) used in the DPG method, we use $\tr{E}$ to avoid a notational conflict with the Fourier transform that is denoted by $\hat{\cdot}$ in the next section.}

\section{Vectorial Envelope DPG Formulation}
\label{sec:envelope}

\paragraph{Mode beat.}
The linear waveguide analysis of the previous section showed how guided waves inside the weakly-guiding step-index fiber fall into a discrete set of LP modes determined by the waveguide characteristics and frequency of the light. Given a specific fiber configuration, each LP mode has its own characteristic wavenumber $\klp$. With the fiber parameters from Table~\ref{tab:fiber-parameters}, the wavenumbers
for the propagating modes are
\begin{equation}
	\{ k_{01}, k_{11}, k_{21}, k_{02} \} \cong
		\{ 8.56833, 8.56630, 8.56380, 8.56322 \}\ \mu\text{m}^{-1} .
	\label{eq:mode-wavenumbers}
\end{equation}
The differences in the modes' propagation constants lead to a distinctive modal interference pattern for each pair of modes, referred to as the \emph{mode beat}.
Let $\Delta \klp := k_{01} - \klp$. Then,
\begin{equation}
	\{ \Delta k_{11}, \Delta k_{21}, \Delta k_{02} \} 
		\cong 
		\{ 2.03, 4.53, 5.11 \}\ \text{mm}^{-1} .
\end{equation}
It follows that the corresponding number of wavelengths per mode beat length between the FM and each HOM respectively is
\begin{equation}
	\frac{ k_{01} }{ \{ \Delta k_{11}, \Delta k_{21}, \Delta k_{02} \} }
	\cong
	\{ 4220.85, 1891.46, 1676.78 \} ,
	\label{eq:mode-beat-length-1}
\end{equation}
which equates to a mode beat length of
\begin{equation}
	\frac{2 \pi}{\Delta \klp} \cong
	\{ 3.10, 1.39, 1.23 \}\ \text{mm} .
	\label{eq:mode-beat-length-2}
\end{equation}

\begin{figure}[htbp]
	\centering
	\begin{subfigure}[b]{1.0\columnwidth}
	\includegraphics[width=\textwidth,trim={60pt 0 10pt 25pt},clip]
	{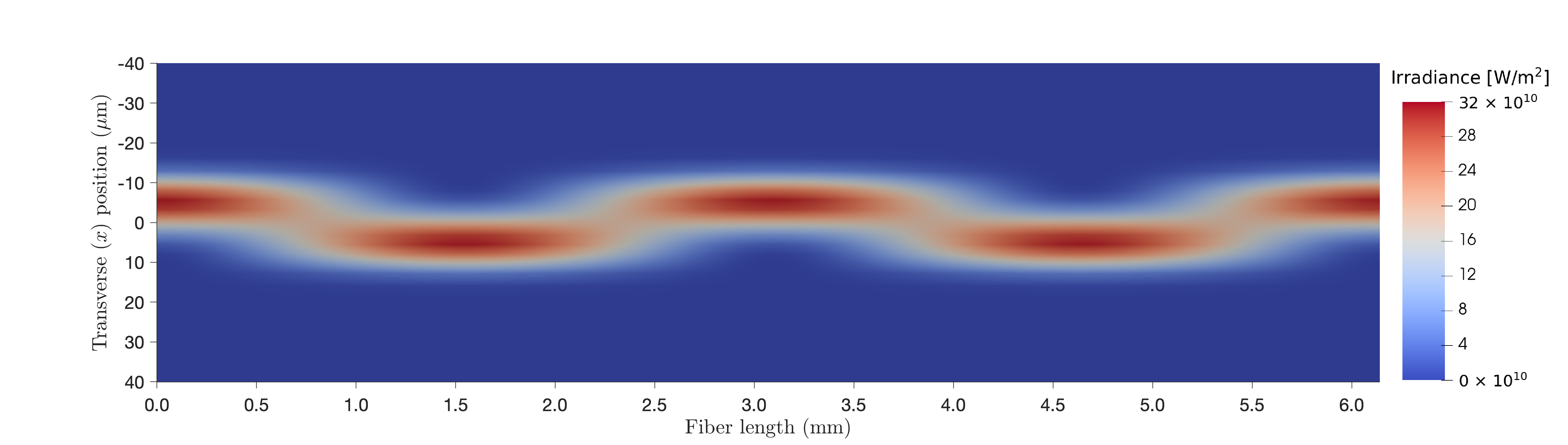}
	\caption*{Mode beat between the $\LP_{01}$ and $\LP_{11}$ modes.}
	\end{subfigure}
	\begin{subfigure}[b]{1.0\columnwidth}
	\includegraphics[width=\textwidth,trim={60pt 0 10pt 25pt},clip]
	{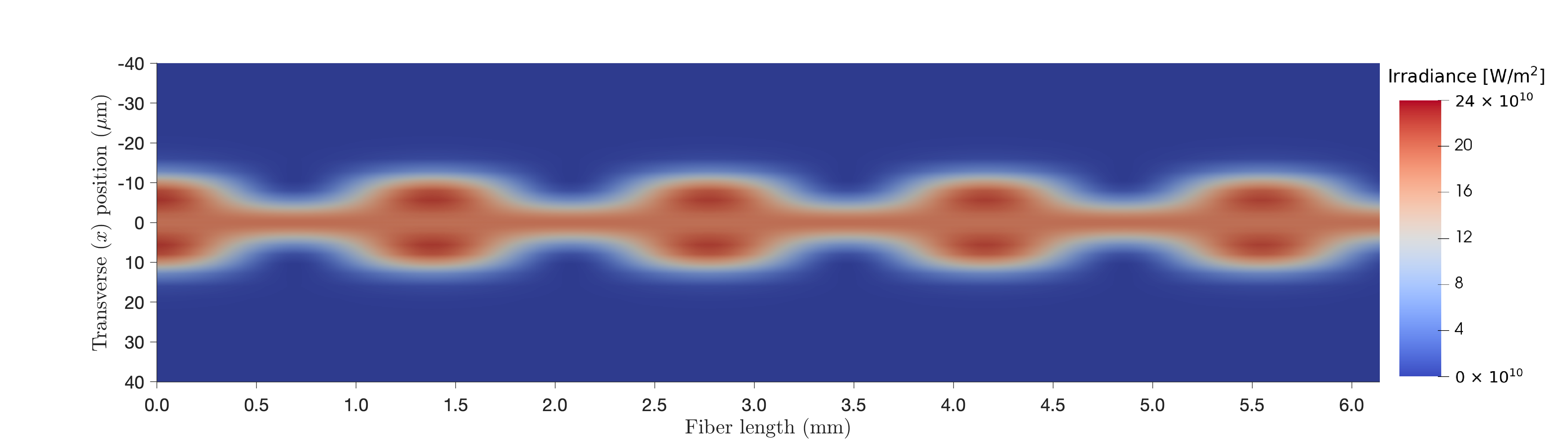}
	\caption*{Mode beat between the $\LP_{01}$ and $\LP_{21}$ modes.}
	\end{subfigure}
	\begin{subfigure}[b]{1.0\columnwidth}
	\includegraphics[width=\textwidth,trim={60pt 0 10pt 25pt},clip]
	{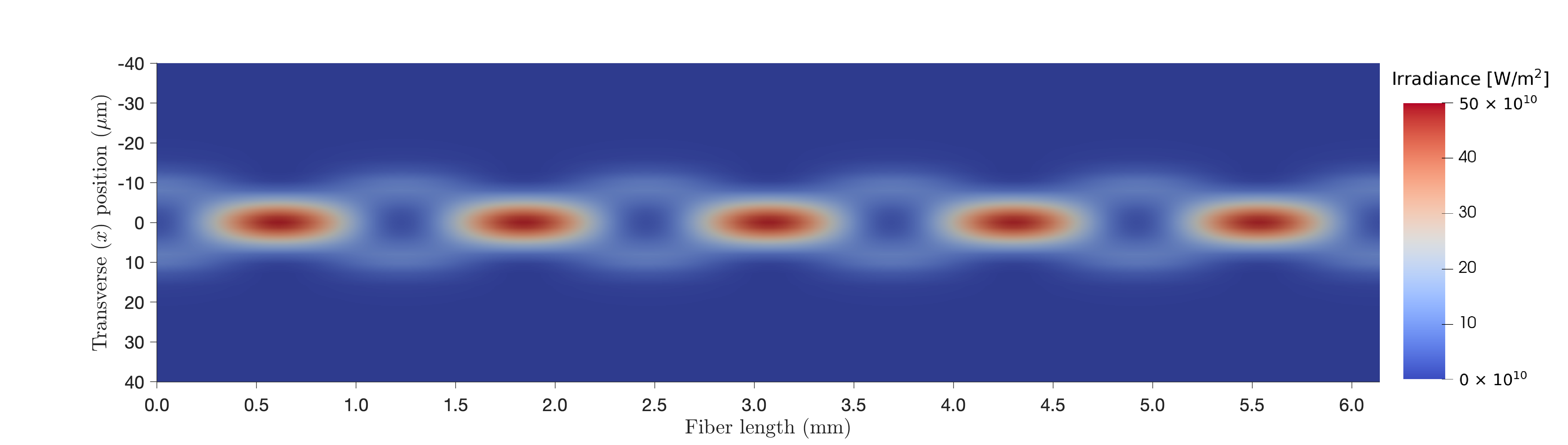}
	\caption*{Mode beat between the $\LP_{01}$ and $\LP_{02}$ modes.}
	\end{subfigure}
	\caption{Irradiance plotted in a longitudinal slice (normal to the $y$-axis) illustrating the mode beat between the FM and the HOMs. The modal interference pattern oscillates at a much longer length scale $(\mc{O}(\text{mm}))$ than the optical wavelength $(\mc{O}(\mu\text{m}))$.}
	\label{fig:mode-beat}
\end{figure}

The mode beat can be nicely visualized by plotting the irradiance, given by $I = | \mf{Re}\{E \times \bar H \} |$ where $\mf{Re} \{ \cdot \}$ denotes the real part of the complex-valued vector, along the longitudinal direction of the fiber, as shown for various combinations of propagating modes in Figure~\fig{mode-beat}. As indicated by \eq{mode-beat-length-1}--\eq{mode-beat-length-2} and the plots in Figure~\fig{mode-beat}, the mode beat occurs on a length scale approximately three orders of magnitudes longer ($\mc{O}$(mm)) than the wavelength ($\mc{O}$($\mu$m)).
  
\paragraph{Envelope ansatz.}
The similarity of the LP mode propagation constants in the weakly-guiding step-index fiber, implying the relatively long mode beat length as compared to the wavelength of the light, motivates the following envelope ansatz for the time-harmonic Maxwell equations:
\begin{align}
	E(x,y,z) &= \ms E(x,y,z) e^{-i \kenv z} ,
	\label{eq:envelope-ansatz-1} \\
	H(x,y,z) &= \ms H(x,y,z) e^{-i \kenv z} ,
	\label{eq:envelope-ansatz-2}
\end{align}
where $\kenv \ge 0$ is the envelope wavenumber, and $\ms E$ and $\ms H$ are the electric and magnetic field envelopes.

At this point, it is important to emphasize the difference between the guided wave assumption \eq{guided-wave-ansatz} that postulates a fixed longitudinal wavenumber $k$ for the guided field envelope $E(x,y)$ and the envelope ansatz \eq{envelope-ansatz-1} that makes no assumptions for the envelope $\ms E(x,y,z)$. In \eq{envelope-ansatz-1}--\eq{envelope-ansatz-2}, the vectorial envelopes $\ms E$ and $\ms H$ are still functions of $(x,y,z)$.

\begin{figure}[!ht]
	\centering
	\begin{subfigure}{\textwidth}
		\includegraphics[width=\textwidth]{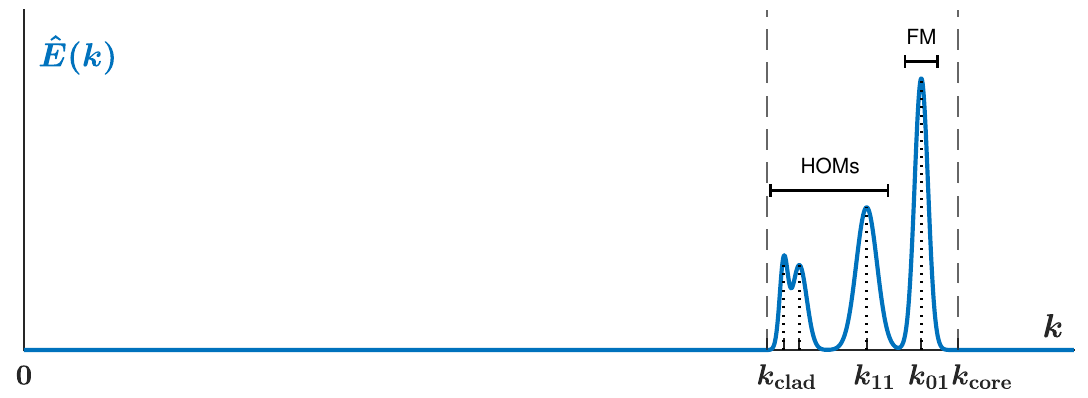}
		\caption{Frequency distribution of the original field $E$ (not to scale)}
		\label{fig:frequency-shift-1}
	\end{subfigure}
	\begin{subfigure}{\textwidth}
		\includegraphics[width=\textwidth]{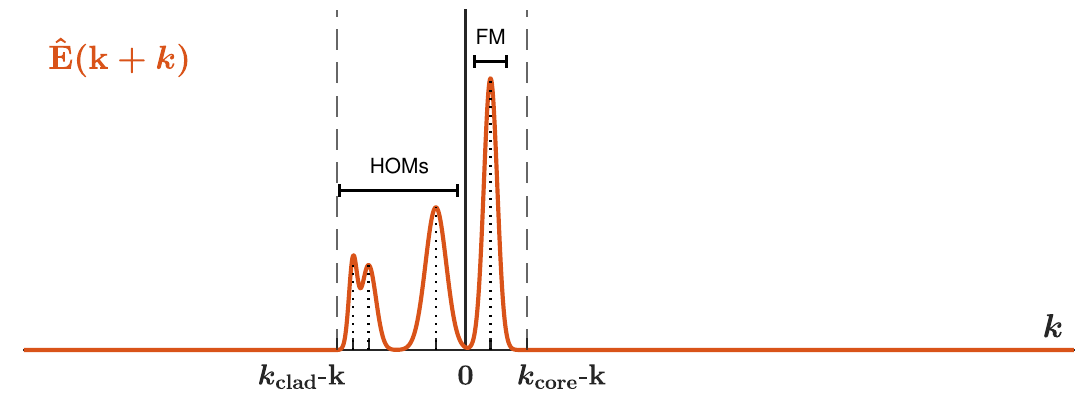}
		\caption{Frequency distribution of the vectorial envelope $\ms E$ (not to scale)}
		\label{fig:frequency-shift-2}
	\end{subfigure}
	\caption{The cost of discretization for resolving the frequencies (spatially in $z$) of the original propagating field $E$, illustrated in (a), is proportional to the maximum frequency $\mc O(\kcore)$; the envelope ansatz shifts the frequencies by the envelope wavenumber $\kenv$, re-centering them nearer the zero frequency ($k=0$), as portrayed in (b), thereby reducing the cost of discretization to the maximum frequency of the envelope which is at most $\mc O(\kcore - \kclad)$. In the weakly-guiding fiber, this frequency shift reduces the maximum frequency by a factor proportional to $\kcore / (\kcore - \kclad) = \mc{O}(1000)$.}
	\label{fig:frequency-shift}
\end{figure}

Effectively, the factor $e^{-i \kenv z}$ in the envelope ansatz re-centers $\ms E$ and $\ms H$ in the frequency domain (via the Fourier transform with respect to $z$) as compared to $E$ and $H$ such that
\begin{equation}
	\mf{F}_{z \to k}\{ E \} =:
	\hat{E}(k) =
	\hat{\ms{E}}(k + \kenv) =
	\mf{F}_{z \to k} \left\{ \ms{E} e^{-i \kenv z} \right\} ,
	\label{eq:frequency-shift}
\end{equation}
where $\mf F$ is the Fourier transform.
This re-centering in the frequency domain is portrayed in Figure~\fig{frequency-shift}: the frequencies of the original field occur within the bandwidth $\kclad < k < \kcore$ (cf.\ Figure~\fig{frequency-shift-1}); the envelope ansatz shifts the frequencies, re-centering them nearer the zero frequency ($k=0$), so that they occur within the bandwidth of $\kclad - \kenv < k < \kcore - \kenv$ (cf.\ Figure~\fig{frequency-shift-2}).

Finite element discretizations are usually dependent on satisfying the Nyquist stability criterion, implying that all propagating wave frequencies must be ``resolved'' to a certain extent in order to have a stable discretization.\footnote{The DPG method provides a stable discretization even if the Nyquist stability criterion is not satisfied; however, accurate solution still requires that the wave is sufficiently resolved~\cite{petrides2021adaptive}.} Typically, this requires a mesh where the element size is on the order of one wavelength or smaller.
By the Nyquist criterion, the cost of discretization is proportional to the highest frequency (or equivalently the total number of wavelengths).\footnote{In fact, the computational cost to obtain an accurate and stable solution may increase even more than linearly with the frequency due to the effect of numerical pollution~\cite{babuska1997pollution}.}
Then, by choosing $\kenv$ in the envelope ansatz \eq{envelope-ansatz-1}--\eq{envelope-ansatz-2} close to the effective wavenumber of the propagating fields $E$ and $H$, such that the respective envelopes $\ms E$ and $\ms H$ are re-centered nearer the zero frequency, the numerical requirements to resolve $\ms E$ and $\ms H$ are significantly reduced.\footnote{Note that if $\kenv$ were chosen not similar to the effective wavenumber, then the envelopes $\ms E$ and $\ms H$ could be just as oscillatory as the fields $E$ and $H$, respectively, rendering the ansatz valid but useless.}

The fact that we are considering laser light implies, by definition, that the optical field has a high degree of spatial (and temporal) coherency, which implies that the field is not broad in frequency space ($k$).
In particular, the propagating field inside the weakly-guiding step-index fiber is comprised of guided LP modes that have similar propagation constants (which is why the mode beat length ($2 \pi / \Delta \klp \sim \mc{O}$(mm)) is several thousand wavelengths). Thus, if the envelope wavenumber $\kenv$ is chosen to be the same as the FM propagation constant, the HOM field envelopes will oscillate in $z$ over length scales of their respective mode beat lengths. The cost of discretizing the envelope may therefore be up to three orders of magnitude lower than a discretization of the original fields.\footnote{In correspondence with the Nyquist criterion, we assume that the cost of discretization increases roughly linearly with the highest frequency of the envelope.}

\begin{remark}{(Slowly varying envelope approximation.)}
It is worth noting that this approach of introducing the envelope ansatz with a propitiously chosen wavenumber $\kenv$ is somewhat similar to the classic \textit{slowly varying envelope approximation} (SVEA) used within the optics field~\cite[{\S}3.3]{shen1984principles}. 
The key idea is that the envelope (field amplitude) is varying at a rate far slower than some pertinent frequency. 
The SVEA uses this fact to neglect specific terms from the optical field governing equations, resulting in a reduced, yet viable, model. 
Indeed, many other models use SVEAs, including the aforementioned coupled-mode-theory TMI models~\cite{naderi2013tmi, goswami2021fiber}. 
Nonetheless, it is important to understand that our methodology does not omit nor neglect terms from the governing equations; it is not an approximation. 
Rather, as will be seen from our results, the Maxwell envelope ansatz~\eqref{eq:envelope-ansatz-1}--\eqref{eq:envelope-ansatz-2} focuses the model around the relevant spatial scales of the mode beat lengths. 
\end{remark}

While the idea of the envelope ansatz \eq{envelope-ansatz-1}--\eq{envelope-ansatz-2} is clearly appealing with regard to the computational efficiency of a step-index fiber model, it remains to demonstrate that this ansatz yields a well-posed formulation of the time-harmonic Maxwell problem \eq{time-harmonic-1}--\eq{time-harmonic-2}. 

The following identities are needed for deriving the envelope equations:
\begin{align}
	\curl (\ms E e^{-i \kenv z}) &= e^{-i \kenv z} (\curl \ms E - i \kenv \ez \times \ms E) ,
	\label{eq:envelope-identity-1} \\
	\ez \times \ms E &= (-\ms E_y, \ms E_x, 0)^T , \\
	\ez \times \ez \times \ms E &= (-\ms E_x, -\ms E_y, 0)^T , \\
	\curl \ez \times \ms E &= (-\ms E_{x,z}, -\ms E_{y,z}, \ms E_{x,x} + \ms E_{y,y})^T , \\
	\ez \times \curl \ms E &= (\ms E_{z,x} - \ms E_{x,z}, \ms E_{z,y} - \ms E_{y,z}, 0)^T ,
	\label{eq:envelope-identity-5}
\end{align}
where $\ms E = (\ms E_x, \ms E_y, \ms E_z)^T$, and $\ez = (0,0,1)^T$ is the unit vector in $z$-direction. The vectorial envelope formulation is obtained by applying the identities \eq{envelope-identity-1}--\eq{envelope-identity-5} to the first-order system \eq{time-harmonic-1}--\eq{time-harmonic-2} with the envelope ansatz \eq{envelope-ansatz-1}--\eq{envelope-ansatz-2}, which yields a modified Maxwell system given by:
\begin{align}
	\curl \ms E - i \kenv \ez \times \ms E &= -i \omega \mu_0 \ms H , 
	\label{eq:envelope-1} \\
	\curl \ms H - i \kenv \ez \times \ms H &= i \omega \eps \ms E .
	\label{eq:envelope-2}
\end{align}

Analogous to \eq{time-harmonic-curl-curl}, a second-order envelope formulation can be derived using the following identity:
\begin{equation}
	\curl \curl (\ms E e^{-i \kenv z}) = e^{-i \kenv z} (\curl \curl \ms E - i \kenv \curl \ez \times \ms E - i \kenv \ez \times \curl \ms E - \kenv^2 \ez \times \ez \times \ms E ) ,
\end{equation}
yielding the curl--curl envelope formulation:
\begin{equation}
	\curl \curl \ms E - i \kenv \curl \ez \times \ms E - i \kenv \ez \times \curl \ms E - \kenv^2 \ez \times \ez \times \ms E - \omega^2 \mu_0 \eps \ms E = 0 .
\end{equation}

In the remainder of the paper, we will limit our discussions to the first-order system \eq{envelope-1}--\eq{envelope-2}.

\paragraph{Ultraweak DPG formulation.}
The discontinuous Petrov--Galerkin (DPG) method \cite{demkowicz2017dpg} yields a pre-asymptotically stable discretization. The built-in stability properties of the method come at the expense of computing \emph{optimal test functions} \cite{demkowicz2011part2}. In practice, this requires the use of discontinuous test spaces leading to the introduction of new unknowns on the mesh skeleton. For a detailed description of the method in general, we refer to \cite{demkowicz2017dpg, gopala2014practical} and references therein; for DPG discretization of time-harmonic Maxwell problems, see \cite{carstensen2016breaking, nagaraj2018raman, henneking2021phd}.

The envelope operator $\ms A : L^2(\Omega) \supset D(\ms A) \rightarrow L^2(\Omega)$ corresponding to the vectorial envelope Maxwell formulation \eq{envelope-1}--\eq{envelope-2} is given by:
\begin{equation}
\ms A \ms u = \left[ \begin{array}{cc}
	-i \omega \eps & (\curl - i \kenv \ez \times) \\
	(\curl - i \kenv \ez \times) & i \omega \mu_0
\end{array} \right]
\left[ \begin{array}{c}
	\ms E \\
	\ms H
\end{array} \right] ,
\label{eq:envelope-operator}
\end{equation}
where $D(\ms A)$ is the domain of the operator incorporating BCs. The corresponding $L^2$ adjoint operator $\ms A^* : L^2(\Omega) \supset D(\ms A^*) \rightarrow L^2(\Omega)$ is given by:
\begin{equation}
\ms A^* \ms v = \left[ \begin{array}{cc}
	-\overline{i \omega \eps} & (\curl + \overline{i \kenv} \ez \times) \\
	(\curl + \overline{i \kenv} \ez \times) & \overline{i \omega \mu_0}
\end{array} \right]
\left[ \begin{array}{c}
	\ms F \\
	\ms G
\end{array} \right] ,
\end{equation}
where the test space $\ms V := D(\ms A^*)$ is equipped with the scaled adjoint graph norm:
\begin{equation}
\begin{split}
	\| \ms v \|_{\ms V}^2 
	&:= \| \ms A^* \ms v \|^2 + \alpha^2 \| \ms v \|^2 \\
	&= \text{the usual terms}\footnotemark + \tc{ORANGE}{\text{the envelope terms}} \\
	&= 
	\| \curl \ms G - \overline{i \omega \eps} \ms F \|^2 + \| \curl \ms F + \overline{i \omega \mu_0} \ms G \|^2
	+ \alpha^2 ( \| \ms F \|^2 + \| \ms G \|^2 ) \\
	&+ 
	\tc{ORANGE}{
	(\curl \ms G - {\overline{i \omega \eps}} \ms F, {\overline{i \kenv}} {\ez \times \ms G})
	+ ({\overline{i \kenv }} {\ez \times \ms G}, \curl \ms G - {\overline{i \omega \eps}} \ms F)
	+ \| {\overline{i \kenv }} {\ez \times \ms G} \|^2 } \\
	&+
	\tc{ORANGE}{
	(\curl \ms F + {\overline{i \omega \mu_0}} \ms G, {\overline{i \kenv }} {\ez \times \ms F})
	+ ({\overline{i \kenv }} {\ez \times \ms F}, \curl \ms F + {\overline{i \omega \mu_0}} \ms G )
	+ \| {\overline{i \kenv }} {\ez \times \ms F} \|^2 } ,
\end{split}
\label{eq:envelope-test-norm}
\end{equation}
\footnotetext{Scaled adjoint graph norm for the standard ultraweak Maxwell problem (without the envelope ansatz).}
and $\alpha \ge 0$ is a scaling coefficient.

The ultraweak formulation of the Maxwell envelope problem \eq{envelope-1}--\eq{envelope-2} is obtained by posing the problem\footnote{Note that $\ms f = 0$ in the waveguide problem in the absence of free charges or impressed currents.}
\begin{equation}
	\left\{ \begin{array}{l}
		\ms u \in D(\ms A) , \\
		\ms A \ms u = \ms f ,
	\end{array} \right.
\end{equation}
in variational form, then integrating by parts both equations and passing all derivatives to the test functions:
\begin{equation}
	\left\{ \begin{array}{l}
		\ms u \in L^2(\Omega) , \\
		(\ms u, \ms A^* \ms v) = (\ms f, \ms v), \quad \ms v \in \ms V .
	\end{array} \right.
\end{equation}

As previously mentioned, in practice the test spaces are broken in the ultraweak DPG formulation to localize the computation of the optimal test functions. This has two important implications: 1)~the test norm \eq{envelope-test-norm} must be localizable, it thus requires a scaling coefficient $\alpha > 0$; and 2)~additional unknowns $\tr{\ms{u}} = (\tr{\ms{E}}, \tr{\ms{H}})^T$ are introduced on the mesh skeleton (see \cite{carstensen2016breaking, henneking2021phd} for details). The broken ultraweak formulation of the vectorial envelope first-order system \eq{envelope-1}--\eq{envelope-2} is given by:
\begin{equation}
	\left\{ \begin{array}{r@{\hskip 3pt}l}
	\multicolumn{2}{l}{\ms E, \ms H \in L^2(\Omega),\ \tr{\ms{E}}, \tr{\ms{H}} \in H^{-1/2}(\tcurl,\Gammah) ,} \\[3pt]
	(\ms E, \hcurl \ms F) + \lb n \times \tr{\ms{E}}, \ms F \rb_{\Gamma_h} + (\ms E, \overline{i \kenv } \ez \times \ms F) + (\ms H, \overline{i \omega \mu_0} \ms F) &= 0, 
\quad \ms F \in \hHcurl , \\[3pt]
	(\ms H, \hcurl \ms G) + \lb n \times \tr{\ms{H}}, \ms G \rb_{\Gamma_h} + (\ms H, \overline{i \kenv } \ez \times \ms G) - (\ms E, \overline{i \omega \eps} \ms G) &= 0, 
\quad \ms G \in \hHcurl ,
	\end{array} \right.
\label{eq:envelope-broken-uw}
\end{equation}
with suitable boundary conditions imposed on the tangential electric field $\tr{\ms{E}}$.

\paragraph{Waveguide stability analysis.}
The stability of the ultraweak formulation for the Maxwell waveguide problem is discussed in \cite{melenk2023waveguide1, demkowicz2024waveguide2}.
The first paper \cite{melenk2023waveguide1} considers a homogeneous waveguide with PEC BCs on $\Gammai$ and $\Gammat$, and a nonlocal Dirichlet-to-Neumann operator that realizes outgoing waves at $\Gammao$.\footnote{The implementation of absorbing BCs is discussed in the next section.} The second paper \cite{demkowicz2024waveguide2} extends the stability result to non-homogeneous waveguides such as the step-index fiber waveguide.\footnote{A third paper considers impedance BCs on one side of the transverse domain for an acoustic waveguide (necessary for modeling loss in bent waveguides) \cite{demkowicz2024waveguide3}.} 

The main result of \cite{melenk2023waveguide1, demkowicz2024waveguide2} is the dependency of the stability constant on the waveguide length~$L$. In particular, the boundedness-below constant $c_b$ of the time-harmonic Maxwell operator $A$ corresponding to \eq{time-harmonic-1}--\eq{time-harmonic-2} is shown to depend linearly upon the waveguide length $L$:
\begin{equation}
	\| A u \| \ge \underbrace{\frac{c_0}{L}}_{=: c_b} \| u \| ,
	\label{eq:bounded-below-constant}
\end{equation}
where $u = (E,H)^T$, and constant $c_0 > 0$ does not depend on $L$.

For the ultraweak formulation with the adjoint graph test norm, the inf--sup constant $\gamma$ depends upon boundedness-below constant $c_b$ and the scaling coefficient $\alpha$ \cite{demkowicz2017dpg}:
\begin{equation}
	\left.
	\begin{array}{ll}
		\| A u \| \geq c_b \| u \| , \, u \in D(A) \\[2pt]
		\| v \|_V^2 = \| A^\ast v \|^2 + \alpha^2 \| v \|^2 
	\end{array}
	\right\}
	\quad \Rightarrow \quad
	\gamma \geq \left[ 1 + \left( \frac{\alpha}{c_b} \right)^2 \right]^{-1/2} .
\end{equation}

The (ideal) DPG method inherits the stability of the continuous problem:\footnote{The practical DPG method only approximates optimal test functions; however, the loss of stability can be controlled by a Fortin operator \cite{gopala2014practical, nagaraj2017fortin}.}$^{,}$\footnote{BA error refers to the best approximation error.}
\begin{equation}
	\underbrace{\| u - u_h \|^2}_{L^2\text{-error}} +
	\underbrace{\| \tr{u} - \tr{u}_h \|_{\tr{U}}^2}_{\text{trace error}} \leq
	\underbrace{ \left[ 1 + \left( \frac{\alpha}{c_b} \right)^2 \right] }_{\text{stability constant}}
	\bigg\{
	\underbrace{\inf_{w_h \in U_h } \| u -w_h \|^2 }_{\text{field BA error }} +
	\underbrace{\inf_{\tr{w}_h \in \tr{U}_h } \| \tr{u} -\tr{w}_h \|_{\tr{U}}^2 }_{\text{trace BA error }}
	\bigg\} ,
\end{equation}
where $\tr{u} = (\tr{E}, \tr{H})^T$, and $\| \cdot \|_{\tr{U}}$ denotes an appropriate minimum energy extension norm for the traces (see \cite{carstensen2016breaking} for details).

As the waveguide length $L$ increases, the loss of stability can thus be countered by inversely proportionally scaling coefficient $\alpha$ in the test norm.\footnote{In practice, the limit for scaling $\alpha$ inversely proportionally is usually determined by round-off errors.} As shown in \cite[Fig.~1]{melenk2023waveguide1}, choosing the right scaling coefficient has a significant impact on the accuracy of the method, which helps to reduce the cost of the discretization in practice.

With the stability of the standard Maxwell waveguide problem established, the stability of the envelope formulation \eq{envelope-1}--\eq{envelope-2} follows by a simple argument \cite[Lem.~2.2]{melenk2023waveguide1}. Recall that the envelope operator $\ms A$, given by \eq{envelope-operator}, is defined as follows:
\begin{equation}
	\ms A \ms u := e^{i \kenv z} A (e^{-i \kenv z} \ms u) .
\end{equation}
Then, using \eq{bounded-below-constant},
\begin{equation}
	\| \ms A \ms u \| =
	\| e^{i \kenv z} A (e^{-i \kenv z} \ms u) \| =
	\| A (e^{-i \kenv z} \ms u) \| \ge
	c_b \| e^{-i \kenv z} \ms u \| =
	c_b \| \ms u \| ,
\end{equation}
implying
\begin{equation}
	\| A u \| \ge c_b \| u \| 
	\Leftrightarrow 
	\| \ms A \ms u \| \ge c_b \| \ms u \| .
\end{equation}
Therefore, the envelope operator $\ms A$ is bounded below if and only if the original operator $A$ is bounded below, and the corresponding boundedness-below constants are identical.

\FloatBarrier

\section{Absorbing Boundary Conditions}
\label{sec:boundary}

The numerical simulation of the fiber waveguide is performed in the truncated computational domain $\Omega$ of finite length $L$. 
At the fiber output $\Gammao$, where the light exits the optical fiber, it is important to avoid unrealistic reflections of the wave that may be caused by imposing (artificial) BCs. 
Absorbing BCs are designed to eliminate or at least mitigate this issue. 

\paragraph{Impedance BC.}
The impedance BC relates the tangential electric and magnetic field unknowns via an impedance relation at $\Gammao$. The enforced impedance relation is of the following form:
\begin{equation}
	H_t = \frac{1}{Z_{\text{imp}}} \ez \times E_t ,
	\label{eq:impedance-BC}
\end{equation}
where $Z_{\text{imp}}$ is the wave impedance which depends on the mode wavenumber and other parameters. 

From the envelope ansatz \eq{envelope-ansatz-1}--\eq{envelope-ansatz-2}, it is easy to see that the impedance relationship \eq{impedance-BC} then also holds for the corresponding envelopes:
\begin{equation}
	\ms H_t = \frac{1}{Z_{\text{imp}}} \ez \times \ms E_t .
	\label{eq:envelope-impedance-BC}
\end{equation}
In other words, the impedance BC is easily adapted to the envelope formulation. In the ultraweak DPG setting \eq{envelope-broken-uw}, the impedance BC \eq{envelope-impedance-BC} is realized using the existing trace unknowns $\tr{\ms E}$ and $\tr{\ms H}$.\footnote{The functional setting on the impedance boundary is not trivial; see \cite{demkowicz2023fem} for details.}

First-order absorbing BCs such as \eq{impedance-BC} are simple to implement and work sufficiently well for many use cases~\cite{engquist1977absorbing}.
In waveguide problems, they can be used to efficiently model single-mode propagation.
However, in the context of multi-mode propagation, such as the simulation of a weakly-guiding step-index fiber with a $V$-number greater than 2.405 (cf.~Table~\tab{LP-cutoff}), or the nonlinear fiber amplifier problem discussed in the next section, it is not sufficient to impose an impedance BC through an impedance relation like \eq{impedance-BC}.

\paragraph{Stretched coordinate PML.}

An effective method for absorbing the propagating wave at the fiber output is the perfectly matched layer (PML). PMLs, first introduced in \cite{berenger1994pml}, are a popular tool in the finite element simulation of wave propagation phenomena in unbounded domains (see \cite{chew1994pml, bramble2007analysis, michler2007pml} and references therein); more recently, stretched coordinate PMLs have been extended to the DPG methodology \cite{astaneh2018pml, nagaraj2018raman, petrides2019phd}.

In the fiber waveguide domain (illustrated in Figure~\fig{fiber-waveguide}), the coordinate stretching is implemented by defining a (uniaxial) complex stretching map:
\begin{equation}
	\bb{R}^3 \ni (x,y,z) \longrightarrow (x,y,\tilde z) \in \bb{C}^3 ,
	\label{eq:stretched-coordinate-pml-1}
\end{equation}
where
\begin{equation}
	\tilde z = \left\{ \begin{array}{ll}
		z &,\ \text{if } 0 < z < l , \\
		z - if (z,\omega) &,\ \text{if } l \leq z < L ,
	\end{array} \right.
	\label{eq:stretched-coordinate-pml-2}
\end{equation}
and $f(z,\omega) > 0$ is the stretching function.\footnote{For example, $f(z, \omega) = (C / \omega) ((z - l) / (L - l))^3$, where constant $C > 0$.}
 
The complex coordinate stretching is designed to cause exponential decay of an outgoing ($k > 0$) guided wave of the form $E e^{-ikz}$. Now suppose the envelope ansatz \eq{envelope-ansatz-1} used wavenumber $\kenv \le \keff$ where the original field $E$ propagated with effective wavenumber $\keff$.\footnote{In the multi-mode step-index fiber waveguide, $\keff \in [ \kclad, \kcore ]$ (cf.~Figure~\fig{frequency-shift-1}).}  Then, applying the complex coordinate stretching to the field envelope $\ms E$ yields:
\begin{equation}
	\ms E \sim e^{-i (\keff - \kenv) \tilde z} = \left\{ \begin{array}{ll}
		e^{-i (\keff - \kenv) z} &,\ \text{if } 0 < z < l , \\
		e^{-i (\keff - \kenv) z} e^{-(\keff - \kenv) f} &,\ \text{if } l \leq z < L .
	\end{array} \right.
	\label{eq:pml-decay}
\end{equation}
Inside the PML region, the exponential decay factor of the field envelope is proportional to $(\keff - \kenv)$ which may be arbitrarily small (or zero if $\keff = \kenv$). If indeed $\keff = \kenv$ (constant envelope), the complex coordinate stretching has no effect. The PML needs to be carefully designed to work effectively with the envelope formulation. We consider two different approaches.

\paragraph{PML envelope formulation 1.}
The first approach follows trivially from \eq{pml-decay}.
If $\tilde k := \keff - \kenv > 0$ is bounded away from zero, then the field envelope decays exponentially, $| \ms E | \sim e^{-\tilde k f}$, inside the PML region.
The stretching function $f$ should be scaled proportional to $\keff / \tilde k$. 
The stretched coordinate PML then works as is, assuming that the PML region is long enough in $z$ (typically several oscillations of the wave).\footnote{Here, the wave propagates as $e^{-i\tilde k z}$.}

Let $J$ denote the (diagonal) Jacobian corresponding to the (uniaxial) complex stretching defined in \eq{stretched-coordinate-pml-1}--\eq{stretched-coordinate-pml-2}, $J = \text{diag}(1, 1, \p \tilde z / \p z)$. Then, after applying the complex-coordinate stretching to the envelope formulation \eq{envelope-1}--\eq{envelope-2} and using Piola transforms for pull-backs to regular Cartesian coordinates (see \cite{astaneh2018pml} for details), the strong form of the pulled-back equations is:
\begin{align}
	|J|^{-1} J \curl \ms E - i \kenv \ez \times J^{-T} \ms E &= -i \omega \mu_0 J^{-T} \ms H , 
	\label{eq:pml-strong-1} \\
	|J|^{-1} J \curl \ms H - i \kenv \ez \times J^{-T} \ms H &= i \omega \eps J^{-T} \ms E .
	\label{eq:pml-strong-2}
\end{align}
Note that due to uniaxial stretching, the Jacobian $J$ has no effect on the rotated components,
\begin{equation}
	\ez \times J^{-T} \ms E = \ez \times \ms E .
\end{equation}
We multiply both equations \eq{pml-strong-1}--\eq{pml-strong-2} with $|J| J^{-1}$; the envelope PML operator and the corresponding $L^2$ adjoint operator are then respectively given by
\begin{equation}
\ms A \ms u
= 
\left[ 
	\begin{array}{cc}
		-\acoef & (\curl - \bcoef \ez \times) \\[2pt]
		(\curl - \bcoef \ez \times) & \ccoef
	\end{array}
\right]
\left[ 
	\begin{array}{c}
		\ms E \\
		\ms H
	\end{array}
\right]
\label{eq:pml-operator}
\end{equation}
and
\begin{equation}
\ms A^* \ms v
= 
\left[ 
	\begin{array}{cc}
		-\overline{\acoef} & (\curl + \overline{\bcoef} \ez \times) \\[2pt]
		(\curl + \overline{\bcoef} \ez \times) & \overline{\ccoef}
	\end{array}
\right]
\left[ 
	\begin{array}{c}
		\ms F \\
		\ms G
	\end{array}
\right] .
\label{eq:pml-adjoint}
\end{equation}

Let ${a}:=\acoef$, ${b}:=\bcoef$, ${c}:=\ccoef$, and $({\ms {F_R}}, {\ms {G_R}})^T := (\ez \times \ms F, \ez \times \ms G)^T$. For the ultraweak formulation, the scaled adjoint graph test norm with the PML, defined analogous to \eq{envelope-test-norm}, is now given by:
\begin{equation}
\begin{split}
	\| \ms v \|_{\ms V}^2 
	&= \text{the usual terms}\footnotemark
	+ \tc{ORANGE}{\text{the envelope terms}} \\
	&= \| \curl \ms G - \bar a \ms F \|^2 
	+ \| \curl \ms F + \bar c \ms G \|^2
	+ \alpha^2 ( \| \ms F \|^2 + \| \ms G \|^2 ) \\
	& + \tc{ORANGE}{
	(\curl \ms G - {\bar a} \ms F, {\bar b} {\ms {G_R}})
	+ ({\bar b} {\ms {G_R}}, \curl \ms G - {\bar a} \ms F)
	+ \| {\bar b} {\ms {G_R}} \|^2 } \\
	& + \tc{ORANGE}{
	(\curl \ms F + {\bar c} \ms G, {\bar b} {\ms {F_R}})
	+ ({\bar b} {\ms {F_R}}, \curl \ms F + {\bar c} \ms G ) 
	+ \| {\bar b} {\ms {F_R}} \|^2 } .
\end{split}
\label{eq:envelope-pml-test-norm}
\end{equation}
\footnotetext{Adjoint graph norm for ultraweak Maxwell with PML without the envelope ansatz.}

\begin{figure}[htb]
	\centering
	\includegraphics[width=\textwidth]{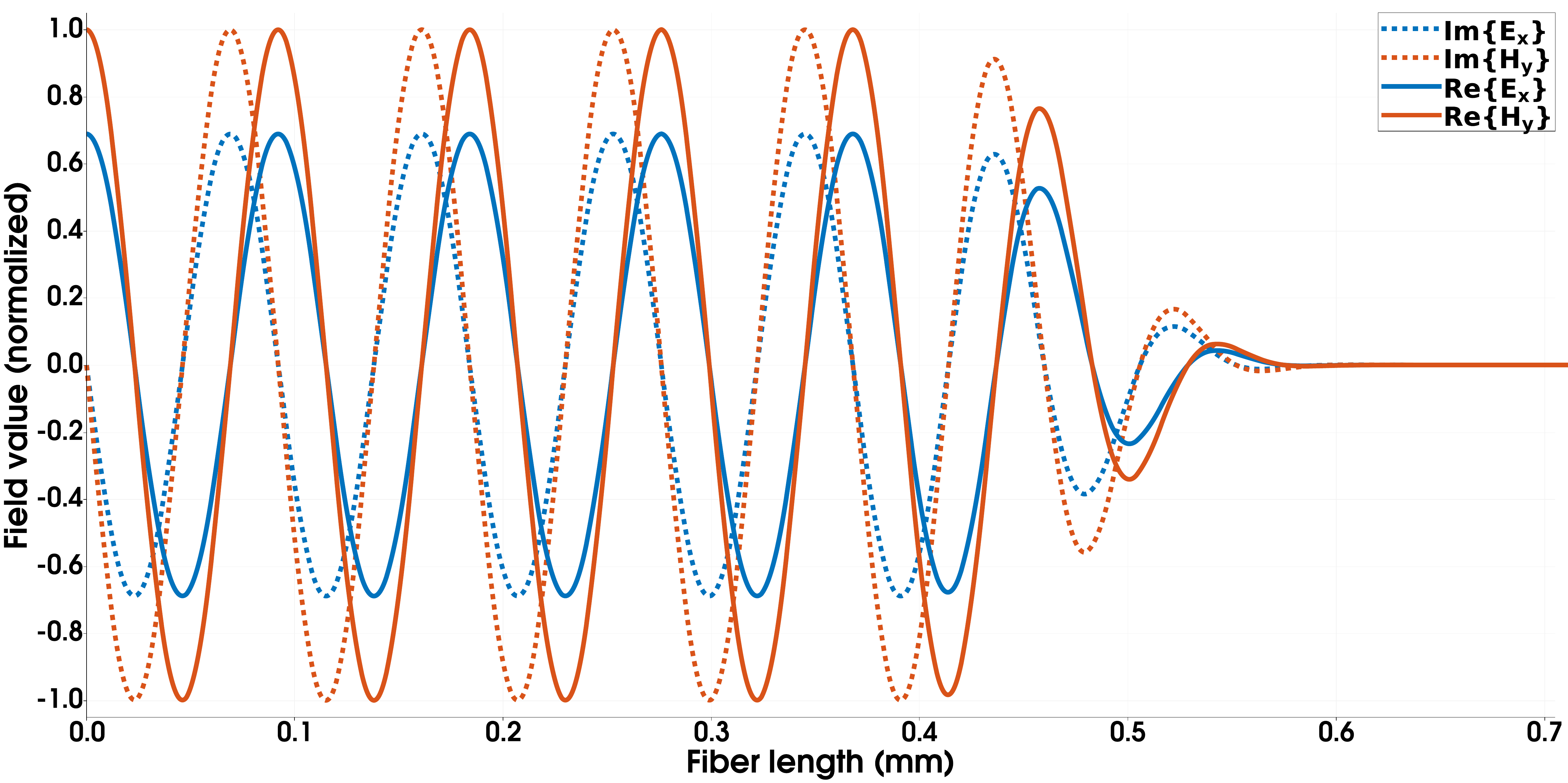}
	\caption{
		Envelope PML formulation 1: single-mode propagation in a step-index fiber waveguide of about $0.7$~mm length ($\sim 1000$ wavelengths). Using an envelope wavenumber of $\kenv = 8.5\ \mu\text{m}^{-1}$, somewhat smaller than the effective wavenumber $\keff = 8.56833\ \mu\text{m}^{-1}$ of the propagating $\LP_{01}$ mode, the envelope is slowly varying along the $z$ direction (ca.~$125$ slower than the original field). Inside the PML region (starting at $z=0.384$~mm and stretching about three envelope beats), the slowly oscillating envelope decays exponentially.}
	\label{fig:envelope-pml-1}
\end{figure}

In practice, it is critical to guarantee a sufficiently large $\tilde k$ for the PML to work. 
This can be ensured by properly choosing $\kenv$ in the envelope ansatz \eq{envelope-ansatz-1}--\eq{envelope-ansatz-2}. 
For example, consider the LP modes in the weakly-guiding step-index fiber. 
For the previous fiber configuration (Table~\tab{fiber-parameters}), the ansatz is chosen based on the mode wavenumbers: 
\begin{equation}
	\kenv < \keff \in \{ k_{01}, k_{11}, k_{21}, k_{02} \} \cong  
		\{ 8.56833, 8.56630, 8.56380, 8.56322 \}\ \mu\text{m}^{-1}.
\end{equation}
Suppose we choose $\kenv = 8.5\ \mu\text{m}^{-1}$; then, $\tilde k = \keff - \kenv \in \{ 68.33, 66.30, 63.80, 63.22 \}\ \text{mm}^{-1}$, yielding an envelope beat length of $\keff / \tilde k \cong \{ 125.40, 125.37, 125.33, 125.32 \}$. That is, the envelope oscillates in $z$ approximately 125 slower than the original field.

Figure~\fig{envelope-pml-1} depicts the transverse electric and magnetic field components of the vectorial envelope for the ($x$-polarized) $\LP_{01}$ mode along the fiber longitudinal ($z$-)axis in a fiber of approximately $1000$ wavelengths. The envelope, oscillating about $125$ more slowly than the original field, decays exponentially once it enters the PML region ($z \ge l = 0.384$~mm). For such a short fiber section, used here for illustrative purposes, the PML region consists of nearly half the domain (ca.~three beats of the envelope).

The advantage of this first PML approach is that it is relatively easy to implement. It uses one consistent envelope ansatz and PML formulation for the entire domain, and the PML domain can be discretized in a similar way to the rest of the domain.\footnote{Using higher-order elements or a finer mesh resolution inside the PML region can still be advantageous for optimizing the PML's performance by reducing numerical errors that can cause artificial reflections of the wave.}
However, using an envelope ansatz with a wavenumber $\kenv$ that is smaller than the effective wavenumber $\keff$, the corresponding electric and magnetic field envelopes are more oscillatory (and thus more expensive to discretize) than strictly necessary.
Recall that the minimum mode beat length in this fiber configuration is more than $1000$ wavelengths (cf.~\eq{mode-beat-length-1}).
Ideally then, the envelope should be roughly three orders of magnitude cheaper to discretize than the original field.

\paragraph{PML envelope formulation 2.}
Our second approach addresses this inefficiency of the first PML formulation by choosing different envelope wavenumbers for the domain of interest $\Omegac := \Omegat \times (0,l)$ and the PML region $\Omegapml := \Omegat \times (l,L)$.
That is, we want to use $\kenv_1 \cong \keff$ in $\Omegac$ to minimize the cost of envelope discretization, and $\kenv_2 < \keff$ in $\Omegapml$ to guarantee exponential decay of the outgoing wave.
By modifying $\kenv$, each domain effectively solves for a different kind of envelope.
In fact, if using $\kenv_2 = 0$ inside the PML region, we would compute the solution to the PML formulation of the standard Maxwell problem.

Varying $\kenv$ locally at the interface of $\Omegac$ and $\Omegapml$, the following interface problems need to be solved:
\begin{equation}
\left\{
\begin{aligned}
	\curl \ms E_1 - i \kenv_1 \ez \times \ms E_1 &= -i \omega \mu_0 \ms H_1 , && z < l, \\
	\curl \ms E_2 - i \kenv_2 \ez \times \ms E_2 &= -i \omega \mu_0 \ms H_2 , && z > l, \\
	\ez \times \ms E_1 e^{-i \kenv_1 z} &= 
	\ez \times \ms E_2 e^{-i \kenv_2 z} , && z = l,
\end{aligned}
\right.
\label{eq:interface-prob-1}
\end{equation}
and
\begin{equation}
\left\{
\begin{aligned}
	\curl \ms H_1 - i \kenv_1 \ez \times \ms H_1 &= i \omega \eps \ms E_1 , && z < l \\
	\curl \ms H_2 - i \kenv_2 \ez \times \ms H_2 &= i \omega \eps \ms E_2 , && z > l \\
	\ez \times \ms H_1 e^{-i \kenv_1 z} & = \ez \times \ms H_2 e^{-i \kenv_2 z} , && z = l.
\end{aligned}
\right.
\label{eq:interface-prob-2}
\end{equation}
Multiplying \eq{interface-prob-1} and \eq{interface-prob-2} with test functions $\ms F_1$, $\ms F_2$, and $\ms G_1$, $\ms G_2$, respectively, and integrating by parts yields:
\begin{align}
	\int_{z < l} \ms E_1 \curl \ms F_1\
	+ \int_{z=l} \ez \times \tr{\ms E}_1 \ms F_1 - 
	\int_{z < l} i \kenv_1 \ez \times \ms E_1 \ms F_1
	&= \int_{z < l} -i \omega \mu_0 \ms H_1 \ms F_1, \\
	\int_{z > l} \ms E_2 \curl \ms F_2\
	- \int_{z=l} \ez \times \tr{\ms{E}}_2 \ms F_2 -
	\int_{z > l} i \kenv_2 \ez \times \ms E_2 \ms F_2
	&= \int_{z > l} -i \omega \mu_0 \ms H_2 \ms F_2 ,
\end{align}
and
\begin{align}
	\int_{z < l} \ms H_1 \curl \ms G_1\
	+ \int_{z=l} \ez \times \tr{\ms H}_1 \ms G_1 - 
	\int_{z < l} i \kenv_1 \ez \times \ms H_1 \ms G_1
	&= \int_{z < l} i \omega \eps \ms E_1 \ms G_1, \\
	\int_{z > l} \ms H_2 \curl \ms G_2\
	- \int_{z=l} \ez \times \tr{\ms H}_2 \ms G_2 -
	\int_{z > l} i \kenv_2 \ez \times \ms H_2 \ms G_2
	&= \int_{z > l} i \omega \eps \ms E_2 \ms G_2 .
\end{align}
The interface conditions can be transferred to the trace unknowns $\tr{\ms E}_1, \tr{\ms E}_2$ and $\tr{\ms H}_1, \tr{\ms H}_2$:
\begin{alignat}{2}
	\label{eq:interf-cond-1}
	\ez \times \tr{\ms E}_1 e^{-i \kenv_1 l} &=
	\ez \times \tr{\ms E}_2 e^{-i \kenv_2 l} &&=:
	\ez \times \tr{\ms E} ,\\
	 \label{eq:interf-cond-2}
	\ez \times \tr{\ms H}_1 e^{-i \kenv_1 l} &=
	\ez \times \tr{\ms H}_2 e^{-i \kenv_2 l} &&=:
	\ez \times \tr{\ms H} .
\end{alignat}
Incorporating the interface conditions \eq{interf-cond-1} and \eq{interf-cond-2} into the formulation, we respectively obtain:
\begin{align}
	\int_{z < l} \ms E_1 \curl \ms F_1 + 
	e^{i \kenv_1 l} \int_{z=l} \ez \times \tr{\ms E} \ms F_1 - 
	\int_{z < l} i \kenv_1 \ez \times \ms E_1 \ms F_1
	&= \int_{z < l} -i \omega \mu_0 \ms H_1 \ms F_1, \\
	\int_{z > l} \ms E_2 \curl \ms F_2 -
	e^{i \kenv_2 l} \int_{z=l} \ez \times \tr{\ms E} \ms F_2 -
	\int_{z > l} i \kenv_2 \ez \times \ms E_2 \ms F_2
	&= \int_{z > l} -i \omega \mu_0 \ms H_2 \ms F_2 ,
\end{align}
and
\begin{align}
	\int_{z < l} \ms H_1 \curl \ms G_1 + 
	e^{i \kenv_1 l} \int_{z=l} \ez \times \tr{\ms H} \ms G_1 - 
	\int_{z < l} i \kenv_1 \ez \times \ms H_1 \ms G_1
	&= \int_{z < l} i \omega \eps \ms E_1 \ms G_1, \\
	\int_{z > l} \ms H_2 \curl \ms G_2 -
	e^{i \kenv_2 l} \int_{z=l} \ez \times \tr{\ms H} \ms G_2 -
	\int_{z > l} i \kenv_2 \ez \times \ms H_2 \ms G_2
	&= \int_{z > l} i \omega \eps \ms E_2 \ms G_2 .
\end{align}
Compared with the standard envelope formulation, the only difference is the presence of the exponential factors in front of the trace terms.
The stretched coordinate PML is then applied analogously to the derivation of \eq{pml-strong-1}--\eq{pml-strong-2}.

\begin{figure}[htb]
	\centering
	\includegraphics[width=\textwidth]{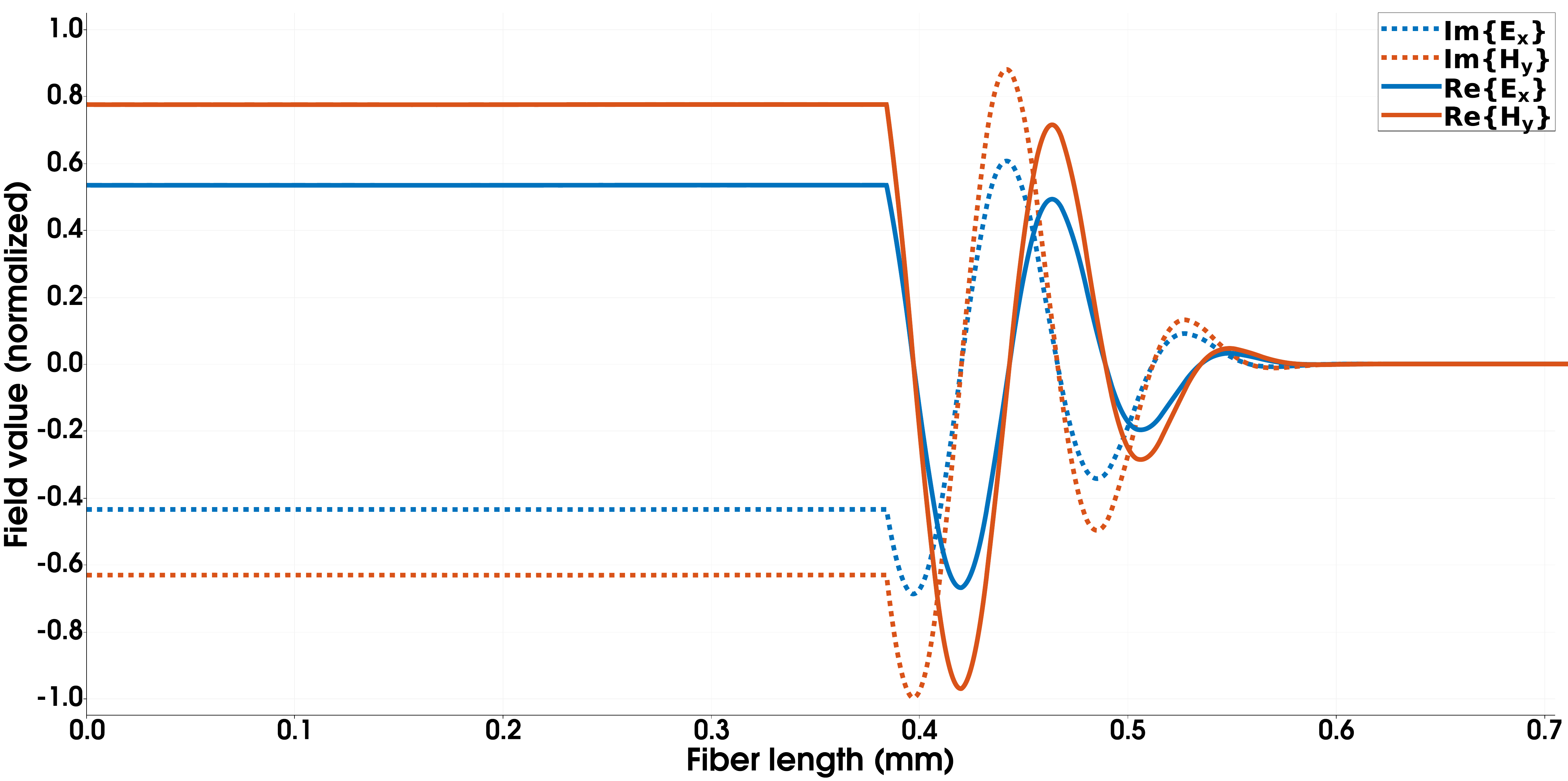}
	\caption{
		Envelope PML formulation 2: single-mode propagation in a step-index fiber waveguide of about $0.7$~mm length ($\sim 1000$ wavelengths). The envelope wavenumber $\kenv_1$ matches exactly the effective wavenumber $\keff = 8.56833\ \mu\text{m}^{-1}$ of the propagating $\LP_{01}$ mode, thus the envelope is constant in the computational domain of interest. Inside the PML region ($z \ge 0.384$~mm), the envelope is computed for $\kenv_2 = 8.5\ \mu\text{m}^{-1}$, enabling exponential decay of the outgoing wave. An interface condition on the traces is embedded in the formulation to accommodate the locally varying wavenumber ansatz.}
	\label{fig:envelope-pml-2}
\end{figure}

Figure~\fig{envelope-pml-2} demonstrates the efficacy of this PML envelope formulation for the single-mode propagation in the step-index fiber. Choosing $\kenv_1 = \keff$ and $\kenv_2 = 8.5\ \mu\text{m}^{-1} < 8.56833\ \mu\text{m}^{-1} = \keff$, the envelope is constant in $\Omegac$ and begins oscillating while exponentially decaying in $\Omegapml$. In the multi-mode propagation case, selecting $\kenv_1 \cong \keff$ implies that the envelope oscillates in $z$ corresponding to the shortest mode beat length.

The approach minimizes the cost of discretization in $\Omegac$. The PML region must still discretize a few envelope beats (with the length depending on the choice of $\kenv_2$) to absorb the wave.
A small drawback of this second approach is that it is perhaps more difficult to implement than the first. 
The changes to the interface terms where $\kenv$ varies locally are easily accomplished; however, the mesh now needs to be adapted locally to the envelope ansatz because of the different wavenumber ansatz.
In particular, the discretization within the PML needs to be fine enough (element size proportional to the wavelength of the envelope in $\Omegapml$), i.e.~potentially much finer than in $\Omegac$ where the element size in $z$ is proportional to the mode beat.

\begin{remark}
Beyond absorbing BCs, the approach of locally varying $\kenv$ can also be useful for other scenarios where the effective wavenumbers of the propagating fields vary along the length of the fiber.
For example, if fiber parameters (such as the core refractive index) vary along the $z$ direction due to the fiber design, heating inside the fiber, or other perturbations, adapting $\kenv$ locally may prove an effective tool for minimizing the cost of discretizing fiber models with the vectorial Maxwell envelope formulation.
\end{remark}

\section{Vectorial Envelope Maxwell Fiber Amplifier Model}
\label{sec:amplifier}

This section discusses briefly how the envelope Maxwell formulation, developed in the previous sections, can be applied to modeling nonlinear laser gain in optical amplifiers.
For this purpose, we consider the 3D vectorial fiber amplifier model developed in \cite{nagaraj2018raman, henneking2021fiber}. The model assumes that the signal laser (wavelength $\lambdas = 1064$~nm) and the pump field (wavelength $\lambdap = 976$~nm) each satisfy the time-harmonic Maxwell equations.

\paragraph{Maxwell fiber amplifier model.}
The amplification of the signal laser occurs through the active gain process in which an active dopant (e.g.~ytterbium) inside the fiber core region absorbs pump light and is subsequently stimulated by the signal laser to emit photons at the signal frequency. This process is modeled as a (first-order) complex perturbation to the refractive index through a weak coupling of the signal and pump Maxwell systems using ion rate equations \cite{pask1995ytterbium}. The strong forms of the weakly-coupled Maxwell systems are given by:
\begin{align}
	\curl \Ej &= -i \omegaj \mu_0 \Hj , 
	\label{eq:coupled-maxwell-1} \\
	\curl \Hj &= i \omegaj (\eps_0 \Ej + \Pj) ,
	\label{eq:coupled-maxwell-2}
\end{align}
where $\Pj$ is the induced electric polarization, which characterizes the light-material interactions, and $\mathrm{j} \in \{\mathrm{s}, \mathrm{p}\}$ refers to either the signal (laser) or pump optical field, respectively. 
The polarization term includes the background polarization,
\begin{equation}
	\Pj^{\text{background}} = \eps_0 (n^2 - 1) \Ej ,
\end{equation}
and the gain polarization,
\begin{equation}
	\Pj^{\text{active gain}} = i \eps_0 \frac{nc}{\omegaj} \gj(\Is, \Ip) \Ej ,
	\label{eq:active-gain}
\end{equation}
where $n$ is the optical fiber's index of refraction, $c$ [m/s] is the speed of causality, $\Ij$ [W/m$^{2}$] are the optical field irradiances, and $\gj$ [1/m] are the active gain functions that model the laser amplification and pump depletion mediated by the ytterbium dopant in the core region (see \cite{henneking2021fiber} for details). 
Thermal polarization is implicitly modeled through a temperature dependence of the  material refractive index. 
Heat deposition inside the fiber amplifier occurs because pump photons have a higher frequency than signal photons; so some energy is lost in the pump-to-signal light conversion by the stimulated gain ($\gj$) process, leading to heat deposition along the fiber. 
This heat deposition can be modeled as a source term in the heat equation~\cite{smith2016mode}: 
\begin{equation}
	Q(\Is, \Ip) = - \Big( \gp(\Is, \Ip) \Ip + \gs(\Is, \Ip) \Is \Big) .
	\label{eq:heat-source}
\end{equation}
For additional details regarding the Maxwell fiber amplifier model, please refer to \cite{nagaraj2018raman, henneking2021fiber, henneking2021phd, henneking2022parallel}. 

The vectorial Maxwell fiber amplifier model can be discretized and solved with the DPG finite element method in the ultraweak setting. 
However, due to the computational expense of resolving the wavelength scale of the light ($\mc{O}(\mu\text{m})$), these computations are limited to fiber lengths of only a few centimeters (tens of thousands of wavelengths) whereas real-length fiber amplifiers are several meters long (millions of wavelengths). 
As will be demonstrated in the next section, applying the envelope formulation to the Maxwell fiber amplifier model helps mitigate these computational impediments such that the full lengths of fiber amplifiers can be simulated. 

\paragraph{Envelope formulation of the Maxwell fiber amplifier model.}

Before applying the vectorial envelope ansatz \eq{envelope-ansatz-1}--\eq{envelope-ansatz-2} to the weakly-coupled Maxwell systems \eq{coupled-maxwell-1}--\eq{coupled-maxwell-2}, we make the following observations.

First, the envelope formulation must be extended to multiple propagating fields.
Each Maxwell system needs to use a different envelope ansatz (i.e.~a different envelope wavenumber $\kenv$) for each distinct optical field (distinct frequency/wavelength):\footnote{In this model, we consider two distinct frequencies---signal ($\omegas$) and pump ($\omegap$)---but the envelope formulation can easily be extended to include additional fields, such as the Stokes wave if modeling Raman gain amplification~\cite{nagaraj2018raman}.}
\begin{align}
	\Ej(x,y,z) &= \ms \Ej(x,y,z) e^{-i \kj z} ,
	\label{eq:coupled-envelope-ansatz-1} \\
	\Hj(x,y,z) &= \ms \Hj(x,y,z) e^{-i \kj z} .
	\label{eq:coupled-envelope-ansatz-2}
\end{align}

Second, for the envelope formulation to be effective, it is necessary that the nonlinear gain and heat source do not depend directly on the highly oscillatory components of the field.
The irradiances computed from the original fields are equivalent to the irradiances computed from the respective vectorial envelopes:
\begin{equation}
	\Ij
	= \left| \mf{Re}\left\{\Ej \times \bar \Hj \right\} \right|
	= \left| \mf{Re}\left\{\ms \Ej e^{-i \kj z} \times \bar{\ms H}_\mathrm{j} e^{i \kj z} \right\} \right|
	= \left| \mf{Re}\left\{\ms \Ej \times \bar{\ms H}_\mathrm{j} \right\} \right|
	= \ms \Ij .
	\label{eq:envelope-irradiance}
\end{equation}
Since the gain function and the heat source each only depend on the fields' irradiances, \eq{envelope-irradiance} implies that:
\begin{align}
	\gj(\Is, \Ip) &= \gj(\ms \Is, \ms \Ip) , 
	\label{eq:envelope-heat-source} \\
	Q(\Is, \Ip) &= Q(\ms \Is, \ms \Ip) .
	\label{eq:envelope-active-gain}
\end{align}
In other words, the gain and heat source can be computed from the slowly varying envelopes $\ms \Ej, \ms \Hj$ instead of the highly oscillatory fields $\Ej, \Hj$.
The weakly-coupled vectorial envelope Maxwell systems of the fiber amplifier model are then given by:
\begin{align}
	\curl \ms \Ej - i\kj \ez \times \ms \Ej &= -i \omegaj \mu_0 \ms \Hj , 
	\label{eq:coupled-envelope-1} \\
	\curl \ms \Hj - i\kj \ez \times \ms \Hj &= i \omegaj \eps_0 
	\left( n^2 \ms \Ej + i \frac{nc}{\omegaj} \gj(\ms \Is, \ms \Ip) \ms \Ej \right) .
	\label{eq:coupled-envelope-2}
\end{align}

The fiber amplifier envelope model \eq{coupled-envelope-1}--\eq{coupled-envelope-2} can be solved via fixed point iteration where the formulations corresponding to signal and pump fields are each discretized with the ultraweak DPG method analogous to the linear vectorial envelope Maxwell formulation \eq{envelope-1}--\eq{envelope-2} presented in Section~\ref{sec:envelope}. The techniques for implementing absorbing BCs presented in Section~\ref{sec:boundary} are also directly applicable to the nonlinear model; however, note that the PML parameters should be optimized for the frequency of the light and the particular envelope ansatz used for each Maxwell system.

\begin{remark}
The fiber amplifier model presented here assumed a co-pumped configuration, where both the signal laser and the pump light are injected at the fiber input $\Gammai$. The model also supports counter-pumped configurations (see \cite{henneking2021phd}) where the pump light is injected at the opposite fiber end $\Gammao$; modeling counter-pumped configurations requires minor adjustments to the envelope model and PML which are not delineated here.
\end{remark}

\section{Numerical Examples: Fiber Amplifier Simulation}
\label{sec:numerical}

The vectorial envelope formulation of the Maxwell fiber amplifier model is implemented in the $hp$3D finite element library\footnote{\url{https://github.com/Oden-EAG/hp3d}}~\cite{henneking2024hp3d, hpUserManual}.
Even with the computational efficiency of the envelope formulation, the Maxwell model is expensive to compute for realistic fiber lengths.
A meter-long optical fiber (millions of optical wavelengths) has thousands of envelope beats that must be resolved accurately by the discretization.
The numerical simulations presented here use high-order ($p \ge 6$) isoparametric elements with curvilinear geometry; the computations are parallelized efficiently over large-scale CPU manycore compute architectures~\cite{hpbook3}.

The numerical examples in this section compute a co-pumped, cladding-pumped amplifier configuration.\footnote{In the cladding-pumped configuration, the pump light is modeled as a plane wave (i.e.~it has no transverse dependence) over the entire core and cladding region, with the pump decay being modeled as an ordinary differential equation for the pump irradiance (as is common in other models~\cite{naderi2013tmi, drake2019equivalent}) instead of solving the computationally expensive pump Maxwell system. The vectorial envelope Maxwell formulation is then only used for modeling the signal laser.}
In each case, the signal laser is seeded into the amplifier with most of its power ($>90\%$) in the $\LP_{01}$ mode and the remaining power in the $\LP_{11}$ mode.

\begin{figure}[htbp]
	\centering
	\begin{subfigure}[b]{1.0\columnwidth}
		\includegraphics[width=1.0\textwidth]{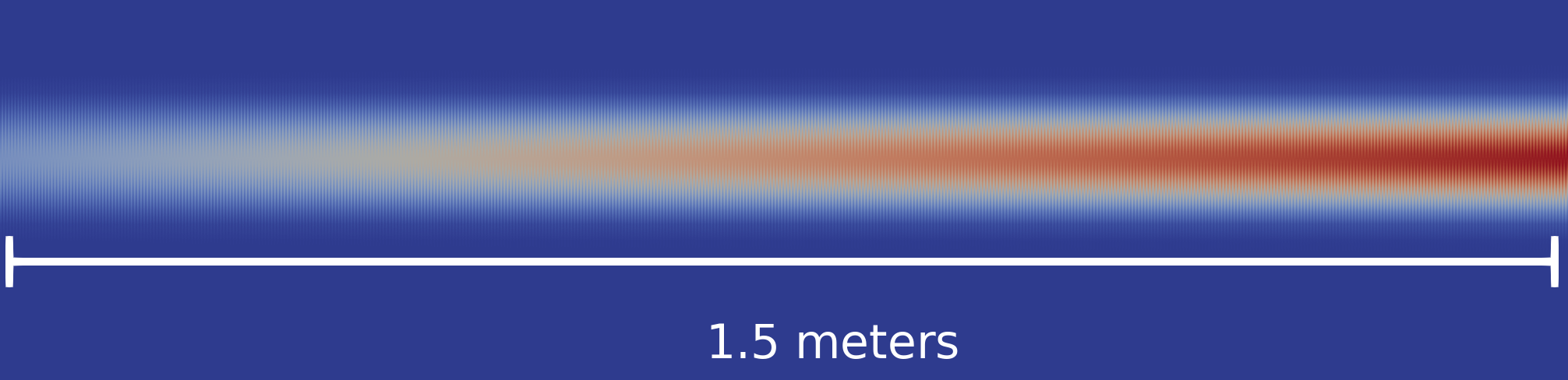}
	\end{subfigure}
	\begin{subfigure}[b]{1.0\columnwidth}
		\includegraphics[width=1.0\textwidth]{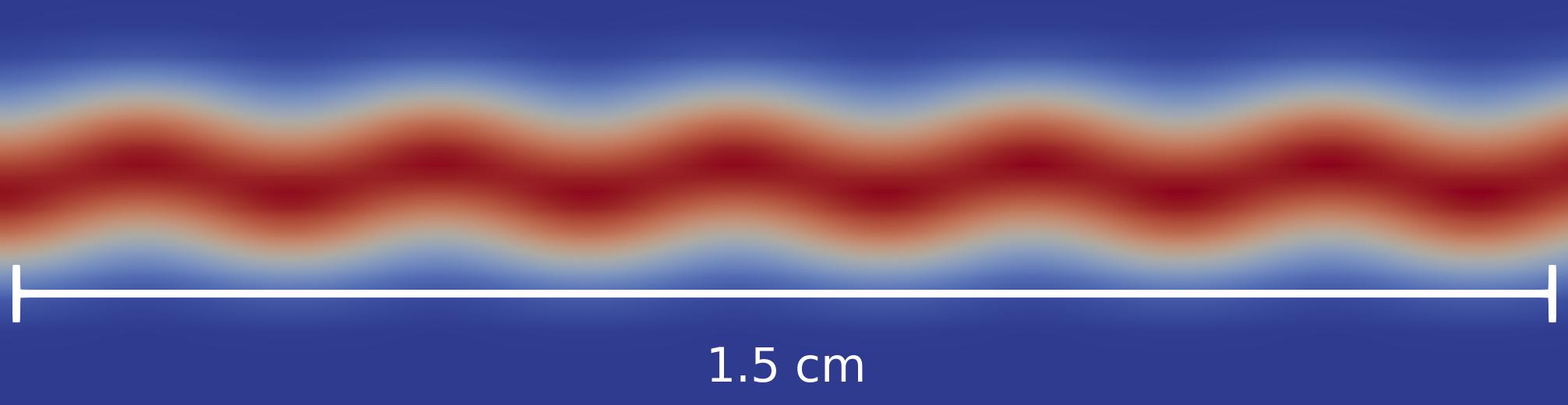}
	\end{subfigure}
	\begin{subfigure}[b]{1.0\columnwidth}
		\includegraphics[width=1.0\textwidth]{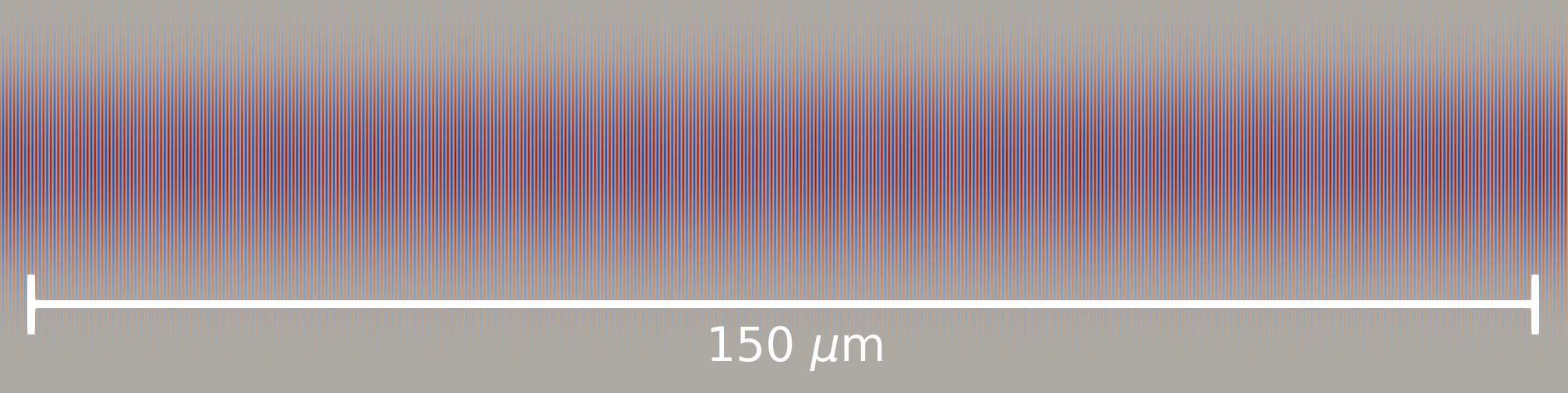}
	\end{subfigure}
	\begin{subfigure}[b]{1.0\columnwidth}
		\includegraphics[width=1.0\textwidth]{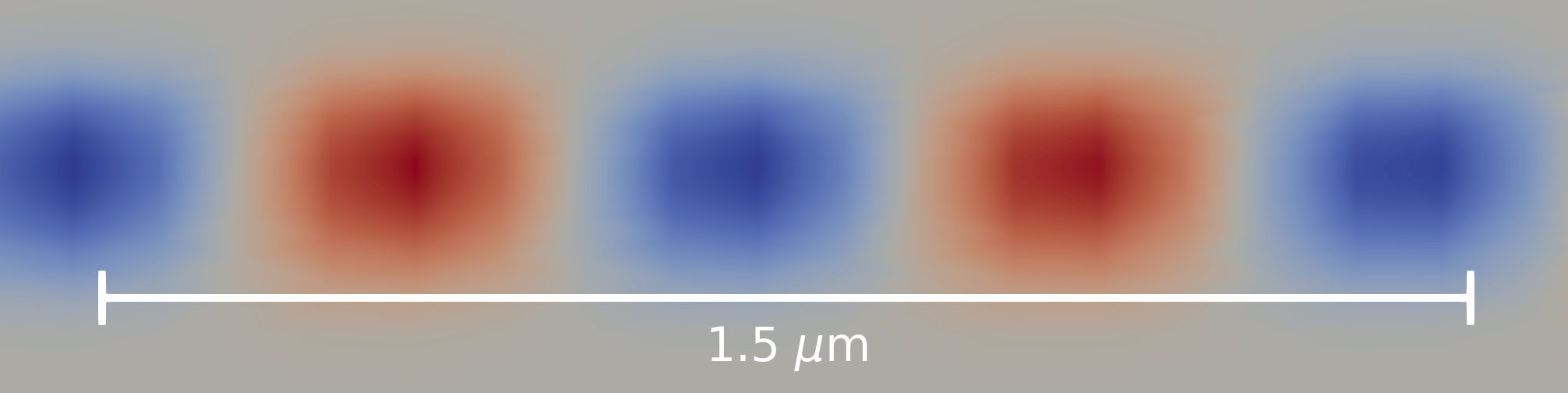}
	\end{subfigure}
	\caption{Numerical simulation of nonlinear laser gain in a 1.6~m long optical waveguide (laser optical wavelength of $\lambdas = 1064$~nm) with the vectorial envelope formulation of the Maxwell fiber amplifier model; each plot shows a longitudinal slice of the fiber waveguide (normal to the $y$-axis).
	 The upper two plots depict the signal irradiance (normalized on a scale of $[0,1]$), which oscillates at the length scale of the mode beat ($\mc O(\text{mm}))$.
	 The lower two plots illustrate the highly oscillatory nature of the electric and magnetic field components (plots show $x$-component of $\mf{Re}\{ \Es \}$ normalized on a scale of $[-1,1]$), visible at the wavelength scale ($\mc O(\mu\text{m}))$.}
	\label{fig:length-scales}
\end{figure}

\paragraph{Nonlinear gain.}
The first example illustrates nonlinear gain inside the fiber amplifier (without heat coupling) for a $1.6$~m long fiber.
Figure~\fig{length-scales} portrays the signal field in various longitudinal slices of the fiber waveguide (normal to the $y$-axis) illustrating various length scales.
Over the full length of the amplifier, the optical power of the signal laser increases significantly through the active gain mechanism, as depicted in the most upper plot of Figure~\fig{length-scales} showing the signal irradiance.
The second plot from the top, also showing signal irradiance, captures the length scale of the mode beat, which is commensurate to the oscillatory behavior of the vectorial envelopes $\ms \Es$ and $\ms \Hs$.
The lower two plots depict the $x$-component of the real part of the electric field $\Es$ that can be retrieved (i.e.~post-processed) from the envelope solution via the ansatz \eq{coupled-envelope-ansatz-1}.
Resolving the laser optical wavelength scale requires a much finer discretization than resolving the envelope, illustrating why computing the 3D vectorial Maxwell model for a full-length fiber amplifier only becomes feasible with the envelope formulation.

\paragraph{TMI simulation.}
Transverse mode instability (TMI) is a thermally-induced nonlinearity. 
TMI is characterized by a sudden reduction of the beam quality above a certain power threshold\footnote{The power threshold depends on the particular fiber laser system, although it is typically anywhere between $100$~W to several kilowatts~\cite{jauregui2020tmi}.} because the various modes of the fiber waveguide start exchanging energy with one another. 
Specifically, there is a transition from a stable beam (below the TMI threshold), where the output signal field is primarily comprised of the FM, to a chaotic energy transfer (above the TMI threshold) between the fiber's guided transverse modes. 
These beam fluctuations happen over millisecond time scales. 
TMI remains a major limitation for the average power scaling of highly coherent fiber laser systems. 
See \cite{jauregui2020tmi, eidam2011experimental} for further discussion of TMI and its origins, characteristics, as well as for known mitigation strategies. 

In order to simulate the onset of TMI for a particular fiber amplifier with a fixed signal seed power, $\Ps(z = 0)$, the launched pump power, $\Pp(z = 0)$, is increased until some figure-of-merit indicates the presence of this nonlinearity. 
We will consider the TMI threshold to have been reached when our metric ($\mtmi$) achieves $5\%$. 
This metric is defined as
\begin{equation}\label{eq:tmi-metric}
	\begin{split}
		\mtmi
		&= \left\langle \sum \Ps^{\text{HOM}}(L,t) / \left( \Ps^{\text{FM}}(L,t) + \sum \Ps^{\text{HOM}}(L,t) \right) \right\rangle_t \\
		&= \left\langle \sum \Ps^{\text{HOM}} (L,t) / \Ps^{\text{total}}(L,t) \right\rangle_t ,
	\end{split}
\end{equation}
where $\lb \cdot \rb_t$ is the time averaging operator. 
Effectively, this figure-of-merit indicates when the signal output from the amplifier, on average, has a significant higher-order mode content, $\sum\Ps^{\text{HOM}}(L,t)$, relative to the total output power, $\Ps^{\text{total}}(L,t)$.

\begin{figure}[htb]
	\centering
	\begin{subfigure}[b]{0.49\columnwidth}
		\includegraphics[width=1.0\textwidth]{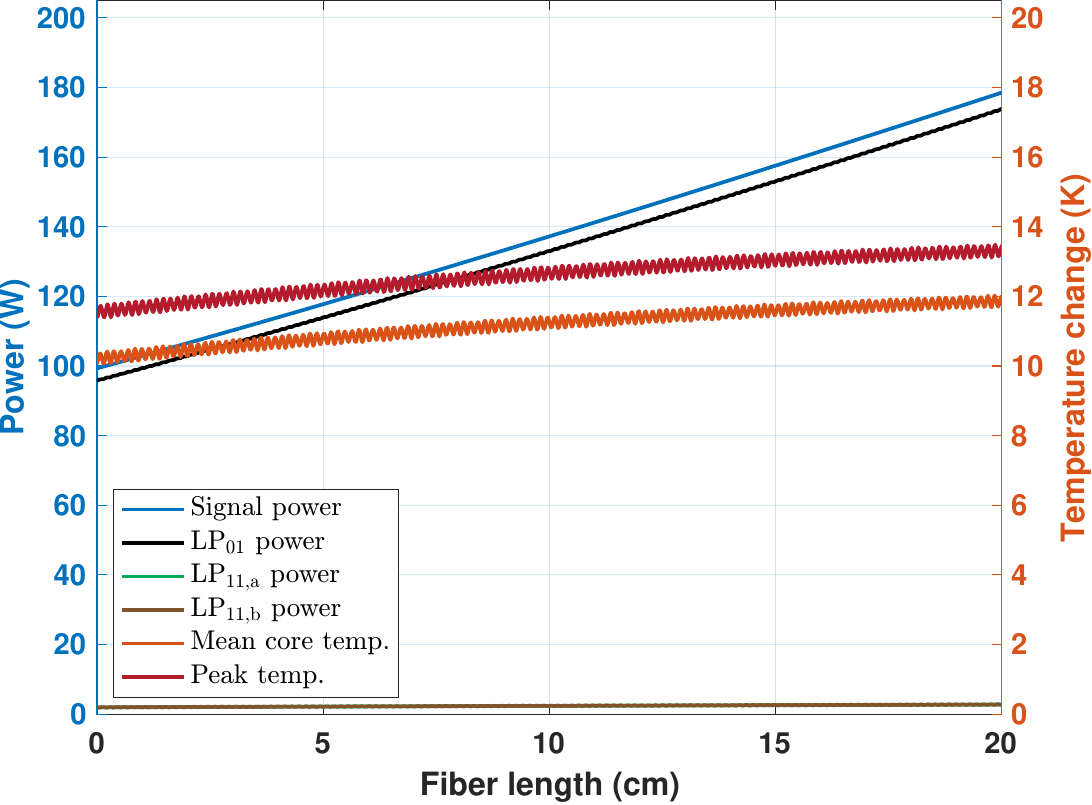}
		\caption{Below TMI threshold, 1~kW pump}
		\label{fig:tmi-stable}
	\end{subfigure}
	\hfill
	\begin{subfigure}[b]{0.49\columnwidth}
		\includegraphics[width=1.0\textwidth]{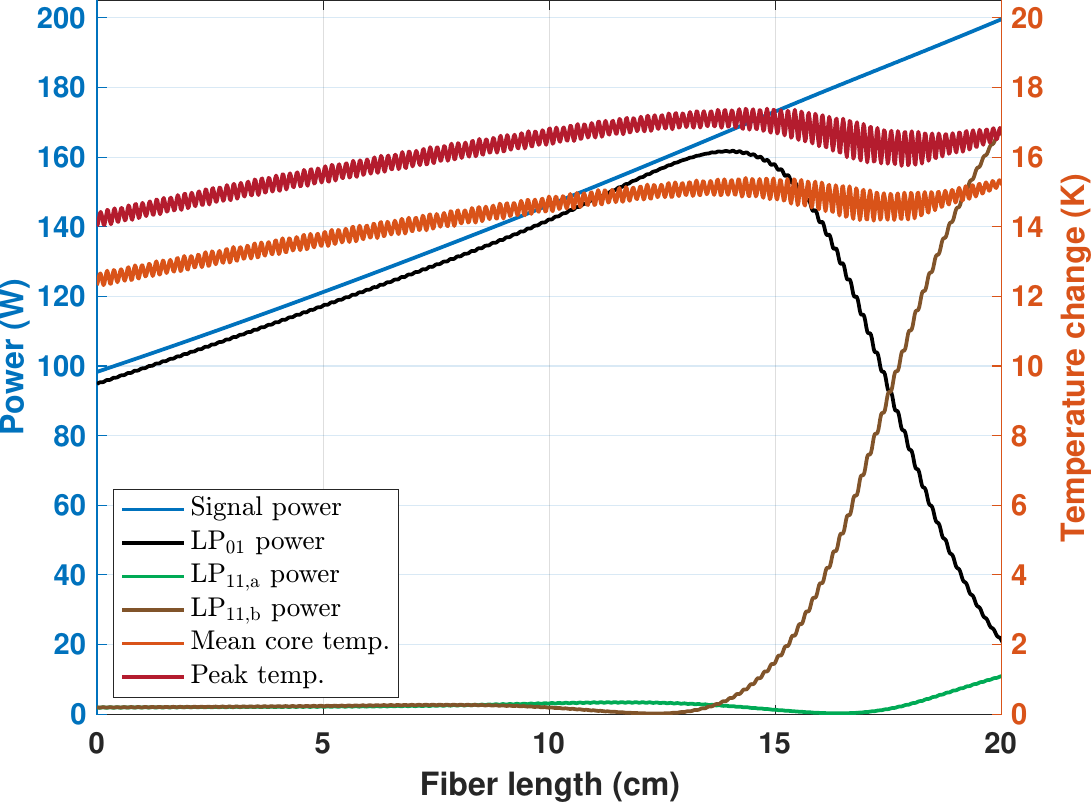}
		\caption{Above TMI threshold, 1.5~kW pump}
		\label{fig:tmi-chaotic}
	\end{subfigure}
	\caption{TMI simulation in a 20~cm long fiber section and a snapshot in time after 25~ms, showing the signal field's mode composition and the fiber's temperature. The (\subref{fig:tmi-stable}) stable regime (below TMI threshold) and the (\subref{fig:tmi-chaotic}) chaotic regime (above TMI threshold) are illustrated.}
	\label{fig:tmi-temperature}
\end{figure}

For TMI simulation, the nonlinear gain model \eq{coupled-envelope-1}--\eq{coupled-envelope-2} is coupled with the transient heat equation using the source term \eq{heat-source}.
Unlike the gain, which occurs over characteristic length scales of $\mc{O}(10\text{ cm})$, capturing the onset of TMI requires simulations that resolve the mode beat lengths, which are $\mc{O}(\text{mm})$.
The mode interference pattern in the signal irradiance translates into a thermally-induced periodic grating of the material refractive index.
This refractive index grating affects each guided transverse mode uniquely, and ultimately induces the energy transfer between the modes, resulting in the TMI.\footnote{The energy transfer between guided transverse modes is linked to a phase shift between the mode interference pattern and the refractive index grating~\cite{jauregui2020tmi}.}
Given the exceedingly large number of wavelengths required to capture sufficiently many mode beats for the TMI to develop, it is extremely challenging to simulate TMI with a vectorial Maxwell model. 
To our best knowledge, the results presented here, using the vectorial envelope formulation, are the first numerical simulations of TMI published for a 3D Maxwell model.

\begin{remark}
As the fiber heats (primarily in the core region), the core numerical aperture increases (known as the \textit{thermal lensing} effect). This tends to separate the mode propagation constants, causing a compression in the mode beat lengths. These mode beat compressions need to be taken into account when discretizing the envelope formulation. Discretizing the envelope in the heated fiber typically requires more discretization points (in $z$) than at ambient temperature. 
\end{remark}

Figure~\fig{tmi-temperature} shows the temperature distribution (as a change from the initial ambient temperature) in a $20$~cm section of the fiber amplifier after $25$~ms of simulation time has elapsed. 
Additionally, the optical power of the signal laser (having been seeded with $100$~W), and its modal composition, are depicted along the fiber longitudinal axis.\footnote{The optical power in each transverse mode is computed from the propagating electromagnetic field using $L^2$~projections onto modes (see~\cite{henneking2021phd}). In general, all of the transverse guided HOMs can carry significant power above the TMI threshold; in this example, nearly all of the energy is in the $\LP_{01}$ and $\LP_{11}$ modes, which is why only those are shown in the plots.}
Figure~\ref{fig:tmi-stable} illustrates the amplifier, with a launched pump power of $1.0$~kW, operating in the stable regime, below the TMI threshold, where nearly all of the signal field's power is found in the $\LP_{01}$ mode (i.e.~the FM).
However, after the launched pump power is increased to $1.5$~kW, as rendered in Figure~\ref{fig:tmi-chaotic}, the amplifier is in the chaotic regime, above the TMI threshold, where energy rapidly transfers between the $\LP_{01}$ and $\LP_{11}$ modes. 
In both plots, the temperature curves exhibit small oscillations associated with the mode beating between these propagating modes. 
Compare from Figure~\ref{fig:tmi-stable} to Figure~\ref{fig:tmi-chaotic} how the amplitude of this refractive index grating is slightly larger at the higher power level, enabling energy exchange between the transverse modes. 
Also, note that the respective amount of power in the modes alters the heat deposition, and thus affects the temperature curve. 

\begin{figure}[htb]
	\centering
	\begin{subfigure}[b]{0.491\columnwidth}
		\includegraphics[width=1.0\textwidth]{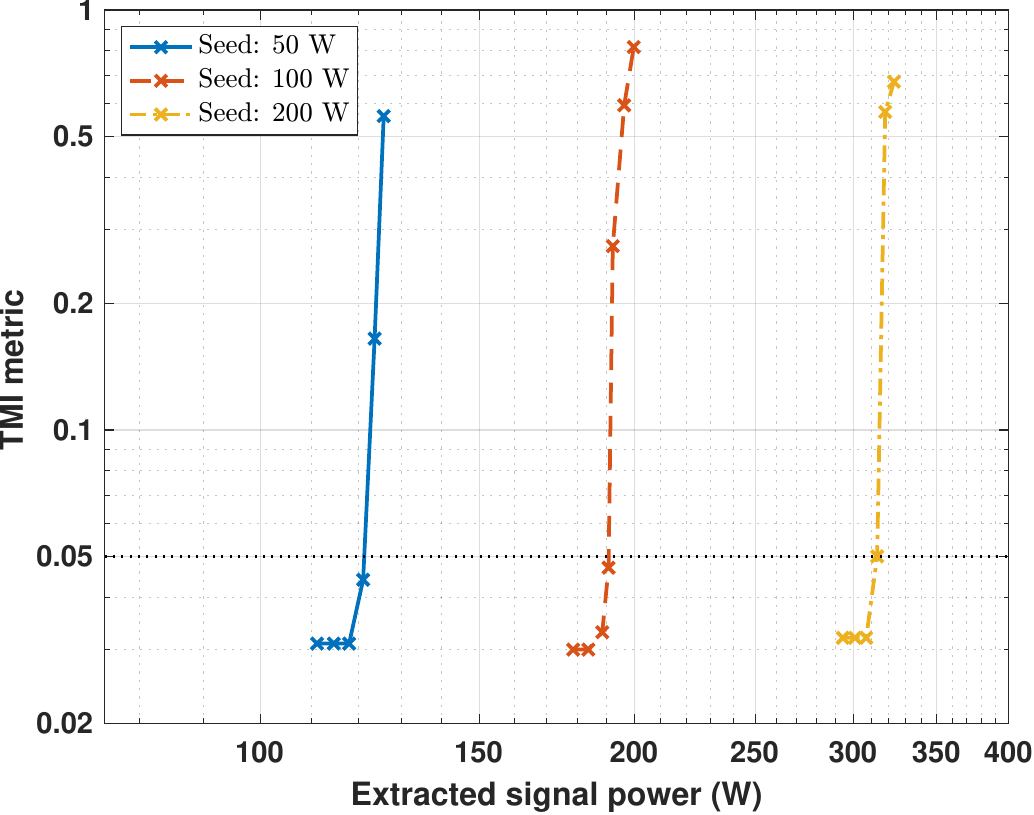}
		\caption{TMI metric}
		\label{fig:tmi-metric}
	\end{subfigure}
	\hfill
	\begin{subfigure}[b]{0.484\columnwidth}
		\includegraphics[width=1.0\textwidth]{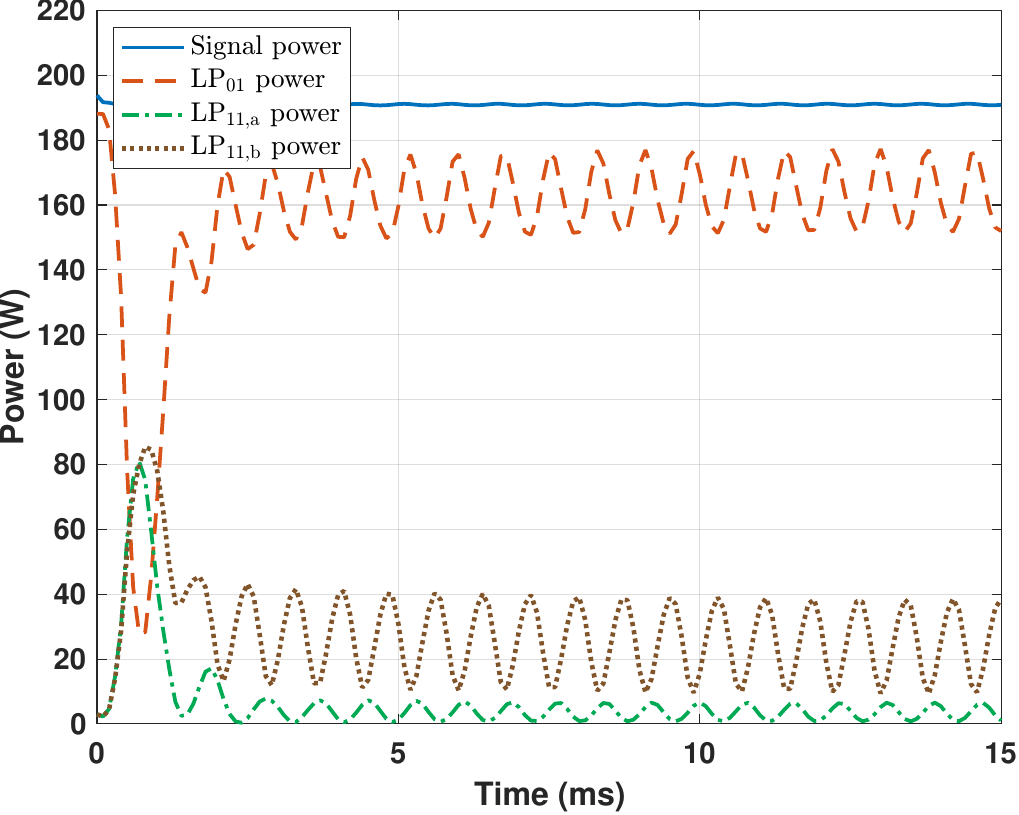}
		\caption{Transition region, 1.29~kW pump}
		\label{fig:tmi-transition}
	\end{subfigure}
	\caption{Related to the TMI simulations of Figure~\ref{fig:tmi-temperature}, the TMI metric~\eqref{eq:tmi-metric} is plotted (\subref{fig:tmi-metric}) as a function of total output signal power ($\Ps(L)$). Note that the launched pump power ($\Pp(0)$) is varied in order to produce each TMI onset curve. These curves represent the onset of the TMI when using three different total signal seed powers ($\Ps(0) \in \{ 50, 100, 200 \}$~W). Near the onset of the TMI threshold, the mode output powers exhibit periodic energy exchanges (\subref{fig:tmi-transition}).}
	\label{fig:tmi}
\end{figure}

Typically, the onset of TMI is observed in the model by tracking the TMI metric~\eqref{eq:tmi-metric} as function of the total extracted signal power from the amplifier. 
Without changing anything else, the output signal power increases as the launched pump power is increased.
Figure~\fig{tmi-metric} shows the onset of the TMI, meaning when $\mtmi = 5\%$ (depicted as a dashed horizontal line) is surpassed, for the three distinct scenarios of $\Ps(0) \in \{ 50, 100, 200 \}$~W.
The increase of the TMI threshold with the signal seed power is a known (nonlinear) dependency (cf.~\cite[Fig.~10b]{jauregui2020tmi}). 
Indeed, there are many design, initial condition, and/or configuration changes that can be made to the fiber amplifier that alter the TMI threshold power, many of which are well-summarized in the TMI review paper~\cite{jauregui2020tmi}.
Below the $\mtmi = 5\%$ threshold, the amplifier is operating in a \textit{stable} regime, as most applications would prefer. 
Once this threshold is surpassed, and up to the initial maxima of the $\mtmi$ function, is the \textit{transition} regime, where modal power exchanges maintain some periodic structure.
With a launched pump power of $1.29$~kW and a $100$~W seed power, Figure~\fig{tmi-transition} illustrates the amplifier operating in the TMI transition regime, where partial, yet mostly periodic, mode power exchanges occur over time.
Afterwards, at even higher launched pump powers, the $\mtmi$ value may actually drop back down to a place that reflects that the output signal power is, on average, equally shared by the various modes of the fiber; this is the \textit{chaotic} regime of the TMI (cf.\ the three regimes~\cite[Fig.~4 of both refs.]{otto2012temporal, naderi2013tmi}).

\section{Conclusions}

In this paper, we proposed a vectorial envelope formulation of the time-harmonic Maxwell equations applicable to computations of electromagnetic waveguides, and to weakly-guiding fiber waveguides in particular.
To accomplish this envelope formulation, one of the more difficult challenges is found in dealing with absorbing BCs.
In particular, we proposed and demonstrated that the PML at the fiber output can be efficiently implemented by leveraging the existing trace unknowns in the ultraweak DPG discretization.
Ultimately, we established that by introducing a sagaciously chosen wavenumber ($\kenv$), our envelope technique re-focuses the standard time-harmonic Maxwell equations to a different characteristic spatial scale, more appropriate for the relevant physics.
In this same spatial direction, the solution of the electromagnetic field inside the fiber waveguide becomes computationally feasible in a much larger domain (in our example, $1000\times$ larger).

To exemplify the success of this method in a complex multiphysics setting, we applied the vectorial envelope formulation to a coupled 3D Maxwell model of an optical fiber laser amplifier. 
By including the active stimulated gain and the appropriate thermal effects in an optical fiber waveguide capable of supporting multiple transverse modes, the model was able to capture the onset of the TMI nonlinearity. 
This was not previously attainable in our past efforts~\cite{henneking2021fiber, henneking2021phd, henneking2022parallel} due to computational constraints imposed by the fact that the standard Maxwell model must resolve the entire fiber down to the wavelength scale ($\mc{O}$($\mu$m)).
Therefore, the computational efficiency of the vectorial envelope methodology makes solution of the 3D Maxwell fiber amplifier model feasible for full-length (meter-long) optical fibers.

\paragraph{Ongoing work.}
The current envelope ansatz is designed to model straight fiber configurations where the fiber longitudinal axis can be aligned with one axis in the Cartesian coordinate system.
In practice, however, fiber amplifiers are typically bent due to spatial and thermal considerations, which impacts nonlinear phenomena including the TMI. 
We therefore aim to extend the current envelope Maxwell formulation to also accommodate circularly coiled fiber configurations. 

The vectorial envelope Maxwell fiber amplifier model is a high-fidelity model in the sense that its derivation requires vastly fewer assumptions than any other models typically used for fiber amplifier modeling.
Given the expense and limitations of obtaining experimental data, this high-fidelity model can be used for validation of computationally efficient lower-fidelity models.
We are currently working toward publishing a detailed comparison of our Maxwell model with a lower-fidelity coupled-mode-theory model for TMI simulation.


\paragraph{Acknowledgements.}
\addcontentsline{toc}{section}{Acknowledgements}
This research was partially supported by AFOSR grants FA9550-19-1-0237 and FA9550-23-1-0103, and NSF Office of Advanced Cyberinfrastructure award 2103524. 
The authors acknowledge the Texas Advanced Computing Center (TACC) at The University of Texas at Austin for providing computational resources under TACC awards DMS22025 (\emph{Frontera} Pathways) and DMS24001 (\emph{Frontera} Leadership Computing Resource Allocations).


\paragraph{Disclaimers.}
\addcontentsline{toc}{section}{Disclaimers}
This article has been approved for public release; distribution unlimited. Public Affairs release approval {\#}AFRL-2024-5730. 
The views expressed in this article are those of the authors and do not necessarily reflect the official policy or position of the Department of the Air Force, the Department of Defense, or the U.S. government. 


\printbibliography[heading=bibintoc]

@string{ACOM="Advances in Computational Mathematics"}

@string{ADVOPT="Advances in Optics and Photonics"}

@string{CAMWA="Computers \& Mathematics with Applications"}

@string{CM="Computational Mechanics"}

@string{CMAM="Computational Methods in Applied Mathematics"}

@string{CMAME="Computer Methods in Applied Mechanics and Engineering"}

@string{FEAD="Finite Elements in Analysis and Design"}

@string{FOCM="Foundations of Computational Mathematics"}

@string{IEEEQ="IEEE Journal of Selected Topics in Quantum Electronics"}

@string{JCP="Journal of Computational Physics"}

@string{JLTEDG="Journal of Lightwave Technology"}

@string{JNA="SIAM Journal on Numerical Analysis"}

@string{JOSAB="Journal of the Optical Society of America"}

@string{JOSS="Journal of Open Source Software"}

@string{MCAMS="Mathematics of Computation"}

@string{MICRO="Microwave and Optical Technology Letters"}

@string{NATPHO="Nature Photonics"}

@string{NM="Numerische Mathematik"}

@string{NMPDE="Numerical Methods for Partial Differential Equations"}

@string{OPTEXP="Optics Express"}

@article{goswami2021fiber,
   title ={Simulations of single-and two-tone {Tm}-doped optical fiber laser amplifiers},
   author = {Goswami, Tathagata and Grosek, Jacob and Gopalakrishnan, Jay},
   journal = OPTEXP,
   volume = {29},
   number = {8},
   pages = {12599--12615},
   year = {2021},
   publisher = {Optical Society of America},
   doi = {10.1364/OE.418095}}

@article{drake2019equivalent,
   author = {Drake, D. and Gopalakrishnan, Jay and Goswami, T. and Grosek, Jacob},
   title ={Simulation of optical fiber amplifier gain using equivalent short fibers},
   journal = CMAME,
   pages = {112698},
   year = {2019},
   publisher = {Elsevier},
   doi = {10.1016/j.cma.2019.112698}}

@article{smith2016mode,
   author = {Smith, A. V. and Smith, J. J.},
   journal = OPTEXP,
   number = {2},
   pages = {975--992},
   publisher = {Optical Society of America},
   title ={Mode instability thresholds for {Tm}-doped fiber amplifiers pumped at 790 nm},
   volume = {24},
   year = {2016},
   doi = {10.1364/OE.24.000975}}

@article{otto2012temporal,
   author = {Otto, H. J. and Stutzki, F. and Jansen, F. and Eidam, T. and Jauregui, Cesar and Limpert, Jens and T\"{u}nnermann, Andreas},
   journal = OPTEXP,
   number = {14},
   pages = {15710--15722},
   publisher = {Optical Society of America},
   title ={Temporal dynamics of mode instabilities in high-power fiber lasers and amplifiers},
   volume = {20},
   year = {2012},
   doi = {10.1364/OE.20.015710}}

@article{jauregui2013fiber,
   author = {Jauregui, Cesar and Limpert, Jens and T\"{u}nnermann, Andreas},
   journal = NATPHO,
   number = {11},
   pages = {861},
   publisher = {Nature Publishing Group},
   title ={High-power fibre lasers},
   volume = {7},
   year = {2013},
   doi = {10.1038/nphoton.2013.273}}

@article{jauregui2020tmi,
   title ={Transverse mode instability},
   author = {Jauregui, Cesar and Stihler, Christoph and Limpert, Jens},
   journal = ADVOPT,
   volume = {12},
   number = {2},
   pages = {429--484},
   year = {2020},
   publisher = {Optical Society of America},
   doi = {10.1364/AOP.385184}}

@book{agrawal,
   author = {Agrawal, Govind P.},
   publisher = {Academic Press},
   title ={Nonlinear fiber optics},
   edition = {5},
   year = {2012},
   isbn = {978-0-12-397023-7}}

@book{griffiths,
   author = {Griffiths, David J.},
   publisher = {Prentice Hall},
   title ={Introduction to electrodynamics},
   edition = {3},
   year = {1999},
   isbn = {013805326X}}

@book{jackson,
   author = {Jackson, John D.},
   publisher = {John Wiley \& Sons},
   title ={Classical electrodynamics},
   edition = {3},
   year = {1999},
   isbn = {0-471-30932-X}}

@book{snyder1983optical,
   title ={Optical waveguide theory},
   author = {Snyder, Allan W. and Love, John D.},
   year = {1983},
   address = {London},
   publisher = {Chapman and Hall Ltd},
   isbn = {0-412-09950-0}}

@book{shen1984principles,
   author = {Shen, Y.-R.},
   publisher = {John Wiley \& Sons},
   title ={The principles of nonlinear optics},
   address = {New York},
   year = {1984},
   isbn = {0-471-88998-9}}

@article{eidam2011experimental,
   author = {Eidam, Tino and Wirth, C. and Jauregui, Cesar and Stutzki, F. and Jansen, F. and Otto, H. J. and Schmidt, O. and Schreiber, T. and Limpert, Jens and T{\"u}nnermann, Andreas},
   journal = OPTEXP,
   number = {14},
   pages = {13218--13224},
   publisher = {Optical Society of America},
   title ={Experimental observations of the threshold-like onset of mode instabilities in high power fiber amplifiers},
   volume = {19},
   year = {2011},
   doi = {10.1364/OE.19.013218}}

@article{gonthier1991bpm,
   author = {Gonthier, F. and H{\'e}nault, A. and Lacroix, S. and Black, R. and Bures, J.},
   journal = JOSAB,
   number = {2},
   pages = {416--421},
   publisher = {Optical Society of America},
   title ={Mode coupling in nonuniform fibers: comparison between coupled-mode theory and finite-difference beam-propagation method simulations},
   volume = {8},
   year = {1991},
   doi = {10.1364/JOSAB.8.000416}}

@article{naderi2013tmi,
   author = {Naderi, S. and Dajani, I. and Madden, T. and Robin, C.},
   journal = OPTEXP,
   number = {13},
   pages = {16111--16129},
   publisher = {Optical Society of America},
   title ={Investigations of modal instabilities in fiber amplifiers through detailed numerical simulations},
   volume = {21},
   year = {2013},
   doi = {10.1364/OE.21.016111}}

@article{pask1995ytterbium,
   author = {Pask, H. M. and Carman, R. J. and Hanna, D. C. and Tropper, A. C. and Mackechnie, C. J. and Barber, P. R. and Dawes, J. M.},
   journal = IEEEQ,
   number = {1},
   pages = {2--13},
   publisher = {IEEE},
   title ={Ytterbium-doped silica fiber lasers: versatile sources for the 1-1.2 $\mu$m region},
   volume = {1},
   year = {1995},
   doi = {10.1109/2944.468377}}

@article{saitoh2001bpm,
   author = {Saitoh, Kunimasa and Koshiba, M.},
   journal = JLTEDG,
   number = {3},
   pages = {405--413},
   publisher = {IEEE},
   title ={Full-vectorial finite element beam propagation method with perfectly matched layers for anisotropic optical waveguides},
   volume = {19},
   year = {2001},
   doi = {10.1109/50.918895}}

@article{ward2013bpm,
   author = {Ward, Benjamin G.},
   journal = OPTEXP,
   number = {10},
   pages = {12053--12067},
   publisher = {Optical Society of America},
   title ={Modeling of transient modal instability in fiber amplifiers},
   volume = {21},
   year = {2013},
   doi = {10.1364/OE.21.012053}}

@article{henneking2024hp3d,
   title = {{$hp$3D}: a scalable {MPI/OpenMP} $hp$-adaptive finite element software library for complex multiphysics applications},
   author = {Henneking, Stefan and Petrides, Socratis and Fuentes, Federico and Badger, Jacob and Demkowicz, Leszek},
   publisher = {The Open Journal},
   journal = JOSS,
   year = {2024},
   volume = {9},
   number = {95},
   pages = {5946},
   doi = {10.21105/joss.05946}}

@article{chakraborty2024hp,
   author = {Chakraborty, Ankit and Henneking, Stefan and Demkowicz, Leszek},
   title = {An anisotropic $hp$-adaptation framework for ultraweak discontinuous {Petrov--Galerkin} formulations},
   journal = CAMWA,
   volume = {167},
   pages = {315--327},
   year = {2024},
   publisher = {Elsevier},
   doi = {10.1016/j.camwa.2024.05.025}}

@article{demkowicz2024waveguide3,
   author = {Demkowicz, Leszek and Gopalakrishnan, Jay and Heuer, Norbert},
   journal = {Oden Institute Report 24-01},
   title = {Stability analysis for acoustic waveguides. {P}art 3: impedance boundary conditions},
   year = {2024}}

@article{demkowicz2024waveguide2,
   author = {Demkowicz, Leszek and Melenk, Jens M. and Badger, Jacob and Henneking, Stefan},
   journal = ACOM,
   title = {Stability analysis for electromagnetic waveguides. {P}art 2: non-homogeneous waveguides},
   volume = {50},
   number = {35},
   year = {2024},
   doi = {10.1007/s10444-024-10130-x}}

@article{melenk2023waveguide1,
   author = {Melenk, Jens M. and Demkowicz, Leszek and Henneking, Stefan},
   title = {Stability analysis for electromagnetic waveguides. {P}art 1: acoustic and homogeneous electromagnetic waveguides},
   year = {2023},
   journal = {arXiv preprint arXiv:2207.12211},
   doi = {10.48550/arXiv.2307.04521}}

@article{badger2023scalable,
   title = {Scalable {DPG} multigrid solver for {Helmholtz} problems: a study on convergence},
   author = {Badger, Jacob and Henneking, Stefan and Petrides, Socratis and Demkowicz, Leszek},
   journal = CAMWA,
   volume = {148},
   pages = {81--92},
   year = {2023},
   publisher = {Elsevier},
   doi = {10.1016/j.camwa.2023.07.006}}

@inproceedings{henneking2022parallel,
   author = {Henneking, Stefan and Grosek, Jacob and Demkowicz, Leszek},
   editor = {Melenk, Jens M. and Perugia, Ilaria and Sch{\"o}berl, Joachim and Schwab, Christoph},
   title = {Parallel simulations of high-power optical fiber amplifiers},
   booktitle = {Spectral and High Order Methods for Partial Differential Equations ICOSAHOM 2020+1},
   year = {2022},
   publisher = {Springer},
   pages = {349--360},
   doi = {10.1007/978-3-031-20432-6}}

@phdthesis{henneking2021phd,
   author = {Henneking, Stefan},
   school = {The University of Texas at Austin},
   title = {A scalable $hp$-adaptive finite element software with applications in fiber optics},
   year = {2021},
   doi = {10.26153/tsw/13716}}

@article{henneking2021pollution,
   title = {A numerical study of the pollution error and {DPG} adaptivity for long waveguide simulations},
   author = {Henneking, Stefan and Demkowicz, Leszek},
   journal = CAMWA,
   publisher = {Elsevier},
   volume = {95},
   pages = {85--100},
   year = {2021},
   doi = {10.1016/j.camwa.2020.03.024}}

@article{henneking2021fiber,
   title = {Model and computational advancements to full vectorial {Maxwell} model for studying fiber amplifiers},
   author = {Henneking, Stefan and Grosek, Jacob and Demkowicz, Leszek},
   journal = CAMWA,
   volume = {85},
   pages = {30--41},
   year = {2021},
   publisher = {Elsevier},
   doi = {10.1016/j.camwa.2021.01.006}}

@article{badger2020fast,
   title = {Sum factorization for fast integration of {DPG} matrices on prismatic elements},
   author = {Badger, Jacob and Henneking, Stefan and Demkowicz, Leszek},
   journal = FEAD,
   publisher = {Elsevier},
   volume = {172},
   pages = {103385},
   year = {2020},
   doi = {10.1016/j.finel.2020.103385}}

@article{hpUserManual,
   author = {Henneking, Stefan and Demkowicz, Leszek},
   title = {{$hp$3D user manual}},
   year = {2022},
   journal = {arXiv preprint arXiv:2207.12211},
   doi = {10.48550/arXiv.2207.12211}}

@book{hpbook3,
   author = {Henneking, Stefan and Demkowicz, Leszek and Petrides, Socratis and Fuentes, Federico and Keith, Brendan and Gatto, Paolo},
   publisher = {In preparation},
   title = {Computing with $hp$ finite elements. {III. Parallel $hp$3D} code},
   year = {2024}}

@book{demkowicz2023fem,
   author = {Demkowicz, Leszek},
   title = {Mathematical theory of finite elements},
   publisher = {Society for Industrial and Applied Mathematics},
   year = {2023},
   doi = {10.1137/1.9781611977738}}

@phdthesis{petrides2019phd,
   author = {Petrides, Socratis},
   School = {The University of Texas at Austin},
   title = {Adaptive multilevel solvers for the discontinuous {Petrov--Galerkin} method with an emphasis on high-frequency wave propagation problems},
   year = {2019},
   doi = {10.26153/tsw/2153}}

@article{petrides2021adaptive,
   title = {An adaptive multigrid solver for {DPG} methods with applications in linear acoustics and electromagnetics},
   author = {Petrides, Socratis and Demkowicz, Leszek},
   journal = CAMWA,
   volume = {87},
   pages = {12-26},
   year = {2021},
   doi = {10.1016/j.camwa.2021.01.017}}

@article{demkowicz2011part2,
   author = {Demkowicz, Leszek and Gopalakrishnan, Jay},
   journal = NMPDE,
   number = {1},
   pages = {70--105},
   publisher = {Wiley Online Library},
   title = {A class of discontinuous {Petrov--Galerkin} methods. {II}. {Optimal} test functions},
   volume = {27},
   year = {2011},
   doi = {10.1002/num.20640}}

@article{engquist1977absorbing,
   title = {Absorbing boundary conditions for numerical simulation of waves},
   author = {Engquist, Bj{\"o}rn and Majda, Andrew},
   journal = MCAMS,
   number = {139},
   pages = {629--651},
   publisher = {American Mathematical Society},
   volume = {31},
   year = {1977},
   doi = {10.1073/pnas.74.5.1765}}

@article{berenger1994pml,
   author = {B{\'e}renger, Jean-Pierre},
   journal = JCP,
   number = {2},
   pages = {185--200},
   publisher = {New York, Academic Press.},
   title = {A perfectly matched layer for the absorption of electromagnetic waves},
   volume = {114},
   year = {1994},
   doi = {10.1006/jcph.1994.1159}}

@article{chew1994pml,
   author = {Chew, W. C. and Weedon, W. H.},
   journal = MICRO,
   number = {13},
   pages = {599--604},
   publisher = {Wiley Online Library},
   title = {A {3D} perfectly matched medium from modified {Maxwell's} equations with stretched coordinates},
   volume = {7},
   year = {1994},
   doi = {10.1002/mop.4650071304}}

@article{bramble2007analysis,
   author = {Bramble, J. and Pasciak, J.},
   journal = MCAMS,
   number = {258},
   pages = {597--614},
   title = {Analysis of a finite {PML} approximation for the three dimensional time-harmonic Maxwell and acoustic scattering problems},
   volume = {76},
   year = {2007},
   doi = {10.1090/S0025-5718-06-01930-2}}

@article{michler2007pml,
   author = {Michler, C. and Demkowicz, Leszek and Kurtz, J. and Pardo, David},
   journal = NMPDE,
   number = {4},
   pages = {832--858},
   publisher = {Wiley Online Library},
   title = {Improving the performance of perfectly matched layers by means of hp-adaptivity},
   volume = {23},
   year = {2007},
   doi = {10.1002/num.20252}}

@article{demkowicz2017dpg,
   author = {Demkowicz, Leszek and Gopalakrishnan, Jay},
   title = {Discontinuous {Petrov--Galerkin} {(DPG)} method},
   journal = {Encyclopedia of Computational Mechanics Second Edition},
   pages = {1--15},
   publisher = {John Wiley \& Sons},
   year = {2017},
   doi = {10.1002/9781119176817.ecm2105}}

@article{melenk2020maxwell,
   title = {Wavenumber-explicit hp-{FEM} analysis for {Maxwell's} equations with transparent boundary conditions},
   author = {Melenk, Jens Markus and Sauter, Stefan A.},
   journal = FOCM,
   pages = {1--117},
   year = {2020},
   publisher = {Springer},
   doi = {10.1007/s10208-020-09452-1}}

@article{babuska1997pollution,
   author = {Babu{\v{s}}ka, Ivo M. and Sauter, Stefan A.},
   journal = JNA,
   number = {6},
   pages = {2392--2423},
   publisher = {SIAM},
   title = {Is the pollution effect of the {FEM} avoidable for the {Helmholtz} equation considering high wave numbers?},
   volume = {34},
   year = {1997},
   doi = {10.1137/S0036142994269186}}

@article{carstensen2016breaking,
   author = {Carstensen, Carsten and Demkowicz, Leszek and Gopalakrishnan, Jay},
   journal = CAMWA,
   number = {3},
   pages = {494-522},
   title = {Breaking spaces and forms for the {DPG} method and applications including {Maxwell} equations},
   volume = {72},
   year = {2016},
   doi = {10.1016/j.camwa.2016.05.004}}

@article{gopala2014practical,
   author = {Gopalakrishnan, Jay and Qiu, W.},
   journal = MCAMS,
   number = {286},
   pages = {537--552},
   title = {An analysis of the practical {DPG} method},
   volume = {83},
   year = {2014},
   doi = {10.1090/S0025-5718-2013-02721-4}}

@article{nagaraj2017fortin,
   author = {Nagaraj, Sriram and Petrides, Socratis and Demkowicz, Leszek},
   journal = CAMWA,
   number = {8},
   pages = {1964--1980},
   publisher = {Elsevier},
   title = {Construction of {DPG Fortin} operators for second order problems},
   volume = {74},
   year = {2017},
   doi = {10.1016/j.camwa.2017.05.030}}

@article{mora2019fast,
   author = {Mora, Jaime David and Demkowicz, Leszek},
   journal = CMAM,
   number = {3},
   pages = {523--555},
   publisher = {De Gruyter},
   title = {Fast integration of {DPG} matrices based on sum factorization for all the energy spaces},
   volume = {19},
   year = {2019},
   doi = {10.1515/cmam-2018-0205}}

@article{nagaraj2018raman,
   author = {Nagaraj, Sriram and Grosek, Jacob and Petrides, Socratis and Demkowicz, Leszek and Mora, Jaime David},
   journal = JCP,
   pages = {100002},
   publisher = {Elsevier},
   title = {A {3D} {DPG} {Maxwell} approach to nonlinear {Raman} gain in fiber laser amplifiers},
   volume = {2},
   year = {2019},
   doi = {10.1016/j.jcpx.2019.100002}}

@article{astaneh2018pml,
   title = {On perfectly matched layers for discontinuous {Petrov--Galerkin} methods},
   author = {Astaneh, Ali V. and Keith, Brendan and Demkowicz, Leszek},
   journal = CM,
   number = {6},
   pages = {1131-1145},
   volume = {63},
   year = {2019},
   doi = {10.1007/s00466-018-1640-3}}


\end{document}